\documentclass[12pt]{article}
\usepackage{comment}
\input epsf
\usepackage{epsfig}
\usepackage{latexsym}
\usepackage{amsfonts}
\usepackage{amsmath}
\usepackage[psamsfonts]{amssymb}

\def\K{{\mathbf{K}}}
\def\H{{\overline{H}}}
\def\sgn{{\mathrm{sgn}}}

\def\R{{\mathbf{R}}}
\def\S{{\mathbf{S}}}
\def\D{{\mathbf{D}}}
\def\P{{\mathbf{P}}}
\def\C{{\mathbf{C}}}

\def\bR{{\mathbf{\overline{R}}}}
\def\P{\mathcal{P}}
\def\Q{\mathcal{Q}}
\def\T{\mathcal{T}}

\newtheorem{thm}{Theorem}[section]
\newtheorem{example}[thm]{Example}
\newtheorem{prop}[thm]{Proposition}
\newtheorem{lemma}[thm]{Lemma}

\newtheorem{cor}[thm]{Corollary}
\newtheorem{rmk}[thm]{Remark}
\newtheorem{df}[thm]{Definition} 

\numberwithin{equation}{section}

\author{Alexandre Eremenko\thanks{Supported by NSF grant DMS-1361836.},
Andrei Gabrielov\thanks{Supported by NSF grant DMS-1161629.}$\, $
and Vitaly Tarasov}
\title{Metrics with four conic singularities and spherical quadrilaterals}
\begin{document}
\maketitle
\begin{abstract} A spherical quadrilateral is a bordered surface
homeomorphic to a closed disk, with four distinguished
boundary points called corners,
equipped with a Riemannian
metric of constant curvature $1$, except at the corners, and such that
the boundary arcs between the corners are geodesic.
We discuss the problem of classification of these quadrilaterals
and perform the classification up to isometry in the case that 
two angles at the corners are multiples of $\pi$.
The problem is equivalent to classification of Heun's equations with real
parameters and
unitary monodromy.

MSC 2010: 30C20,34M03.

Keywords: surfaces of positive curvature, conic singularities, Heun equation,
Schwarz equation, accessory parameters, conformal mapping, circular polygons
\end{abstract}

\noindent
\section{Introduction}\label{sec:Intro}

Let $S$ be a compact Riemann surface, and $a_0,\ldots,a_{n-1}$ a finite set
of points on $S$. Let us consider a conformal Riemannian metric on $S$
of constant curvature
$K\in\{0,1,-1\}$ with conic singularities at the points $a_j$.
This means that in a local conformal coordinate $z$ the length element of the
metric is given by the formula $ds=\rho(z)|dz|$, where
$\rho$ is a solution of the differential equation
\begin{equation}\label{0}
\Delta\log\rho+K\rho^2=0\quad{\mbox{in}}
\quad S\backslash \{ a_0,\ldots,a_{n-1} \},
\end{equation}
and $\rho(z)\sim |z|^{\alpha_j-1}$, for the local coordinate $z$ which
is equal to $0$ at $a_j$. Here $\alpha_j>0,$ and $2\pi\alpha_j$ is
the total angle around the singularity $a_j$.

Alternatively, for every point in $S$, there exists a local coordinate
$z$ for which
$$ds=\frac{2\alpha|z|^{\alpha-1}|dz|}{1+K|z|^2},$$
where $\alpha>0$.
At the singular points $a_j$ we have $\alpha=\alpha_j$ while at
all other points, $\alpha=1$.

In 1890, G\"ottingen Mathematical society proposed the study of
equation (\ref{0}) as a competition topic, probably by
suggestion of H.A.~Schwarz \cite{SG}.
E.~Picard wrote several papers on the
subject, \cite{P1,P2,P}, see also \cite[Chap. 4]{P3}.
When $K<0$, the topic is closely related to uniformization of orbifolds
\cite{SG,For,Po}, where one is interested in the
angles $2\pi\alpha_j=2\pi/m_j$ with positive integer $m_j$.
The case $K\leq 0$ is quite well understood, but very little is known on
the case $K>0$.

McOwen \cite{MO} and
Troyanov \cite{T} studied the general question of existence and uniqueness
of such metrics with prescribed $a_j$, $\alpha_j$ and $K$.
Troyanov also considered
the case of non-constant curvature $K$.
One necessary condition that one has to impose on these data follows
from the Gauss--Bonnet theorem:
the quantity
\begin{equation}\label{curv}
\chi(S)+\sum_{j=0}^{n-1}(\alpha_j-1)\quad\mbox{has the same sign as}\quad K.
\end{equation}
Here $\chi$ is the Euler characteristic.
Indeed, this quantity multiplied by
$2\pi$ is equal to the integral curvature of the smooth part of the surface.

It follows from the results of Picard, McOwen
and Troyanov, that for $K\in\{0, -1\}$,
condition (\ref{curv})
is also sufficient for the existence of the metric with conic singularities
at arbitrary points $a_j$ and angles $2\pi\alpha_j$. The metric with
given $a_j$ and $\alpha_j$ is unique
when $K=-1$, and unique up to a constant multiple when $K=0$.

In the case of positive curvature, the results are much less complete.
The result of Troyanov that applies to $K=1$ is the following:
\vspace{.1in}

{\em Let $S$ be a compact Riemann surface, $a_0,\ldots,a_{n-1}$ points
on $S$, and
\newline
$\alpha_0,\ldots,\alpha_{n-1}$ po\-si\-ti\-ve num\-bers
sa\-tis\-fy\-ing
\begin{equation}\label{troy}
0<\chi(S)+\sum_{j=1}^n(\alpha_j-1)<2\min\{1,\min_{0\leq j\leq n-1}\alpha_j\}.
\end{equation}
Then there exists a conformal metric of positive curvature $1$ on $S$
with conic singularities at $a_j$ and angles $2\pi\alpha_j$.}
\vspace{.1in}

F.~Luo and G.~Tian \cite{LT} proved that if the condition
$0<\alpha_j<1$ is satisfied, then (\ref{troy}) is necessary and sufficient,
and the metric with given $a_j$ and $\alpha_j$ is unique.

In general, the right hand side inequality in (\ref{troy})
is not a necessary condition, and the metric may not be unique.

In this paper, we only
consider the simplest case when $S$ is the sphere, so $\chi(S)=2$.

The problem of description and classification of conformal metrics
of curvature $1$ with conic singularities on the sphere has
applications to the study of certain surfaces of constant mean curvature
\cite{FKKRUY,D2,D}, and to several other questions of
geometry and physics \cite{HR,T}.

The so-called ``real case'' is interesting and important.
Suppose that all singularities belong to a circle on the sphere $S$,
and we only
consider the metrics which are symmetric with respect to this circle.
Then the circle splits $S$ into two symmetric disks. Each of them
is a {\em spherical polygon}, (surface)
for which we state a formal definition:

\begin{df}\label{n-gon}{\rm {\em  A spherical $n$-gon} is a closed disk
with $n$ distinguished boundary points $a_j$ called the corners, equipped with
a conformal Riemannian metric of constant curvature
$1$ everywhere except the corners, and such that the sides (boundary arcs
between the corners) are geodesic. The metric has conic singularities at
the corners.}
\end{df}

\begin{example}\label{flat-n-gons}
{\rm Let us consider {\em flat} $n$-gons, which are defined
similarly.
The necessary and sufficient
condition for the existence of a flat $n$-gon with angles $\pi\alpha_j$
is given by the Gauss--Bonnet theorem which in this case says that
$$\sum\alpha_j=n-2,$$
and the polygon with given angles and prescribed corners\footnote{By this we
mean that the images of the corners under a conformal map onto
a disk are prescribed.}
is unique up to a scaling factor. The simplest proof of these facts
is the Schwarz--Christoffel formula.
Thus our subject can be considered as a generalization of the
Schwarz--Christoffel formula to the case of positive curvature.}
\end{example}
\vspace{.1in}

In \cite{UY,E}, all possibilities
for spherical triangles are completely described,
see also \cite{FKKRUY} where a minor error in \cite[Theorem 2]{E}
was corrected. In the case of triangles,
the metric is uniquely determined by the angles when none of the $\alpha_j$ is an integer.

The case when all $\alpha_j$ are integers, and $n$ is arbitrary,
is also well understood.
In this case, the line element of the metric has the global representation
$$ds=\frac{2|f'||dz|}{1+|f|^2},$$
where $f$ is a rational function. The singular points $a_j$ are the critical
points of $f$, $\alpha_j-1$ is the multiplicity
of the critical point $a_j$, and $\alpha_j$ is the local degree
of $f$ at $a_j$.

Thus the problem with all integer $\alpha_j$ is equivalent to describing
rational functions with prescribed critical points \cite{G,S,EG0,EG,EG2,EGSV}.

Almost nothing is known in the case when
some of the $\alpha_j$ are not integers, the number of singularities
is greater than $3$, and the right-hand side
inequality in (\ref{troy}) is violated.

In this paper we begin investigation of the case
$n=4$, with the emphasis on the real case.

In the real case,
we may assume without loss of generality that the circle is the
real line $\R\cup\{\infty\}$.
The real line splits the sphere $S$ into two symmetric spherical
$n$-gons with the real
corners $a_j$ and the angles
$\pi\alpha_j$ at the corners.

Thus we arrive at
the problem of classification of spherical quadrilaterals (surfaces).
We study this problem using two different methods.
One of them is the geometric method
of F. Klein \cite{Klein} who classified spherical triangles.
Klein classified not only  spherical triangles with
geodesic sides, but also circular triangles, whose sides have
constant geodesic curvature (that is, locally they are arcs of circles).
A modern paper which uses Klein's approach to triangles is \cite{Y}.

Classification of triangles permitted Klein
to obtain exact relations for the numbers of
zeros of hypergeometric functions on
the intervals between the singular points on the real line.
Van Vleck \cite{VV} extended this approach of Klein,
using the same geometric
method, and obtained exact inequalities for the numbers of zeros
of hypergeometric functions in the upper and lower half-planes.
Hurwitz \cite{H} re-proved these results with a different, analytic
method.

We hope that our results can be used to obtain information about
solutions of Heun's equation,
in the same way as Klein obtained information
about solutions of the hypergeometric equation.

Klein's classification of triangles was partially extended to
arbitrary circular quadrilaterals (not necessary geodesic ones) in
the work of Sch\"onflies \cite{Sch,Sch2} and Ihlenburg \cite{I,I2}.
They considered certain geometric reduction process of cutting
a circular quadrilateral into simpler ones. Then they classified the
irreducible quadrilaterals up to  conformal automorphisms of the sphere.
Thus they obtained an algorithm which
permits to construct all circular quadrilaterals. Using
this algorithm, Ihlenburg derived relations between the angles and sides
of a circular quadrilateral. However this algorithm falls short of
a complete classification. In particular, a quadrilateral with prescribed
angles and sides is not unique. Moreover, it seems difficult to single out
geodesic quadrilaterals in the construction of Sch\"onflies and Ihlenburg.

We use somewhat different approach which consists in associating
to every spherical geodesic quadrilateral a combinatorial object
which we call a net, thus reducing the classification to combinatorics.
Then we solve this combinatorics problem and obtain a classification
of spherical quadrilaterals up to isometry. Our approach can be also applied
to general (non-geodesic) circular quadrilaterals.

The boundary of a spherical polygon is a closed curve on the sphere,
consisting of geodesic pieces. The angles of such a curve at the corners
make sense only modulo $2\pi$. All possible sequences of angles
have been described by Biswas \cite{B}, see also \cite{REU}.
These inequalities on the angles give necessary but not sufficient conditions
that the angles of a spherical polygon (surface) must satisfy.

Our second method is a direct study of Heun's equation.
Classification of spherical quadrilaterals can be stated in
terms of a special eigenvalue problem for this equation.
This method leads to complete results when the eigenvalue problem
can be solved algebraically.

The contents of the paper is the following.
In section \ref{DE-connection},
we recall the connection of the problem with Heun's equation, and recall
the results on Heun's equation related to our problem.
In section \ref{all-integers} we describe
what is known when all $\alpha_j$ are integers,
with a special emphasis on the ``real case'' when all $a_j$ belong
to the real line.

In section \ref{cases-leq-3} we
begin the study of the case when two
of the $\alpha_j$ are integers and two others are not.
(If three of the $\alpha_j$ are integers then all four must be integers.)
Complete classification for this case is obtained in the remaining sections.
Sections \ref{cond-angles}-\ref{counting-real} are based on a direct
study of the eigenvalue problem for Heun's equation.
The results are illustrated with numerical examples in section \ref{sec:examples}.
Sections \ref{intro-to-part-2}-\ref{chains} are based on
geometric and combinatorial methods.
The cases when three or four of the $\alpha_j$ are not integers are
postponed to a forthcoming paper.

\section{Connection with linear differential\newline equations}\label{DE-connection}

Let $(S,ds)$ be the Riemann sphere equipped
with a metric with conic singularities.
Every smooth point of $S$ has a neighborhood which is isometric to a region on
the standard unit sphere $\S$; let $f$ be such an isometry.
Then $f$ has an analytic continuation along every path
in $S\backslash\{ a_0,\ldots,a_{n-1}\}$, and we obtain a multi-valued function
which is called the {\em developing map}. The monodromy of
$f$ consists of orientation-preserving isometries (rotations) of
$\S$, so the Schwarzian derivative
\begin{equation}
\label{schwarz}
F(z):=\frac{f^{\prime\prime\prime}}{f^\prime}-\frac{3}{2}\left(
\frac{f^{\prime\prime}}{f^\prime}\right)^2
\end{equation}
is a single valued function.

Developing map is completely characterized by the properties that
it has an analytic continuation along any curve in
$S\backslash\{ a_1,\ldots,a_n\}$,
has asymptotics $\sim c(z-a_j)^{\alpha_j},\; z\to a_j,\; c\neq 0,$ and has
$PSU(2)=SO(3)$ monodromy. It is possible that two such maps with the same
$a_j$ and $\alpha_j$ are related by
post-composition with a fractional-linear transformation.
The metrics arising from such maps will be called {\em equivalent}.
Following \cite{FKKRUY}, we say that the metric is reducible
if its monodromy group is commutative (which is equivalent to
all monodromy transformations having a common fixed point).
In the case of irreducible metrics, each equivalence class
contains only one metric. For reducible metrics, the equivalence
class is a
one-parameter family when the monodromy is non-trivial and
a two-parameter-parametric family when monodromy is trivial.

The asymptotic
behavior of $f$ at the singular
points $a_j$ implies that the only singularities of $F$ on
the sphere are double poles, so $F$ is a rational function,
and we obtain the Schwarz differential equation (\ref{schwarz}) for $f$.

It is well-known that the general solution
of the Schwarz differential equation
is a ratio of two linearly independent solutions of
the linear differential
equation
\begin{equation}
\label{linear}
y^{\prime\prime}+Py'+Qy=0,\quad f=y_1/y_0,
\end{equation}
where $$F=-P'-P^2/2+2Q.$$
For example one can take $P=0$, then $Q=F/2$.
Another convenient choice is to make all poles but one of $P$ and $Q$
simple.
When $n=3$, equation (\ref{linear})
is equivalent to the hypergeometric equation,
and when $n=4$ to Heun's equation \cite{R}.

The singular points $a_j$ of the metric
are the singular points of the equation
(\ref{linear}). These singular points are regular,
and to each point correspond
two exponents $\alpha_j^\prime>\alpha_j^{\prime\prime}$,
so that $\alpha_j=\alpha_j^\prime-\alpha_j^{\prime\prime}$.
If $\alpha_j$ is an integer for some $j$,
we have an additional condition
of the absence of logarithms in the formal solution
of (\ref{linear}) near $a_j$.

It is easy to write down the general form of a Fuchsian equation with
prescribed
singularities and prescribed exponents at the singularities.
After a normalization, $n-3$ parameters remain,
the so-called accessory parameters. To obtain a conformal metric
of curvature $1$, one has to choose these accessory parameters
in such a way that the monodromy group of the equation is conjugate to a subgroup
of $PSU(2)$.

By a fractional-linear change of the independent variable,
one can place one sin\-gu\-lar point at $\infty$.
Then, making
chan\-ges of the va\-riable $y(z)\mapsto y(z)(z-a_j)^{\beta_j}$,
one can assume that the smaller expo\-nent at each finite singular point
is $0$, see \cite{R}.
For the case of four singularities $a_0,\dots,a_3$, where $a_3=\infty$,
we thus obtain Heun's equation in the standard form
\begin{equation}\label{heun}
y^{\prime\prime}+\left(\sum_{j=0}^2\frac{1-\alpha_j}{z-a_j}\right)y^{\prime}+
\frac{Az-\lambda}{(z-a_0)(z-a_1)(z-a_2)}y=0,
\end{equation}
where
\begin{equation}\label{cond}
A=\alpha^\prime\alpha^{\prime\prime},\quad
\sum_{j=0}^2\alpha_j+\alpha^\prime+\alpha^{\prime\prime}=2.
\end{equation}
Here the exponents at the singular points are described
by the Riemann symbol
$$P\left\{\begin{array}{ccccc}a_0&a_1&a_2&\infty&\\
0&0&0&\alpha^{\prime\prime}&;z\\
\alpha_0&\alpha_1&\alpha_2&\alpha^{\prime}&
\end{array}\right\}.$$
The first line lists the singularities, the second the smaller
exponents, and the third the larger exponents.
So the angle at infinity is $\alpha_3=\alpha^\prime-\alpha^{\prime\prime}$.
The accessory parameter is $\lambda$.

Solving the second equation (\ref{cond}) together with
$\alpha'-\alpha^{\prime\prime}=\alpha_3$, we obtain
\begin{equation}
\label{alphaprime}
\alpha'=\frac{1}{2}(2+\alpha_3-\alpha_0-\alpha_1-\alpha_2)
\end{equation}
and
$$\alpha^{\prime\prime}=\frac{1}{2}(2-\alpha_3-\alpha_0-\alpha_1-\alpha_2).
$$
The question of the existence of a spherical quadrilateral
with given corners $a_0$, $a_1$, $a_2$, $\infty$
and given angles $\pi\alpha_j,\; 0\leq j\leq 3$,
is equivalent to the following:
{\em when one can choose real $\lambda$ so that the monodromy
group of Heun's equation (\ref{heun})
is conjugate to a subgroup of $PSU(2)$ ?}

The necessary condition (\ref{curv}) can be restated for the equation
(\ref{heun})
as
\begin{equation}
\label{alphapp}
\alpha^{\prime\prime}<0.
\end{equation}
We also have
$$A=\alpha^\prime\alpha^{\prime\prime}=
\alpha^{\prime\prime}(\alpha_3+\alpha^{\prime\prime})$$
by a simple computation.

One can write (\ref{heun}) in several other forms. Assuming that all $a_j$
are real, we have the Sturm-Liouville form:
\begin{equation}\label{sl}
\frac{d}{dx}\left(\prod_{j=0}^2|x-a_j|^{1-\alpha_j}\, y^\prime\right)+
\frac{(Ax-\lambda)\, \sgn\left(\prod_{j=0}^2(x-a_j)\right)}{\prod_{j=0}^2|x-a_j|^{\alpha_j}}y=0.
\end{equation}
Sometimes the Schr\"odinger form is more convenient:
\begin{equation}\label{schro}
y^{\prime\prime}
-\left(\lambda+\frac{1}{4}\sum_{k=0}^3
\frac{\alpha_k^2-1}{x-a_k}
\prod_{j\neq k}(a_k-a_j)\right)\frac{y}{\prod_{j=0}^3(x-a_j)}=0,
\end{equation}
where all four singularities are in the finite part of the plane.
The exponents in the Schr\"odinger form are $(1\pm\alpha_j)/2.$
The potential in (\ref{schro}) is  $F/2$, where
$F$ is the Schwarzian (\ref{schwarz}).

A question similar to our problem
was investigated in \cite{Klein2,Hilb1,Hilb2,Sm,Sm2}:
{\em when can one choose the accessory parameter so that the monodromy group
of Heun's equation preserves a circle?}
All these authors consider the problem under the assumption
\begin{equation}
\label{assumption}
0\leq\alpha_j<1,\quad\mbox{for}\quad 0\leq j\leq 3.
\end{equation}
The most comprehensive treatment of this problem is in Smirnov's thesis
\cite{Sm}. Smirnov proved that for all sets of data satisfying
(\ref{assumption}), there exists a sequence of values of
the accessory parameter $\lambda=\lambda_k,\; k=0,\pm1,\pm2,\ldots$ such
that the monodromy group of the equation has an invariant circle.
Each of the two opposite sides of the corresponding quadrilateral covers
a circle
$|k|$ times, and the other two sides are proper subsets of their
corresponding circles.

The problem of choosing the accessory
parameter so that the monodromy group is conjugate to a subgroup in $PSU(2)$ is discussed in \cite{D}.
However all results of that
paper are also proved only under the assumption (\ref{assumption}).

Assumption (\ref{assumption}) seems to be essential for the methods
of Klein \cite{Klein2}, Hilb, Smirnov and Dorfmeister.

\section{The case of all integers corners}\label{all-integers}

If all $\alpha_j$ are integers, the developing map $f$ is a rational
function, and the metric of curvature $1$
with conic singularities can be globally described
as the pull-back of the spherical metric via $f$,
that is
$$ds=\frac{2|f'||dz|}{1+|f|^2}.$$
The singular points $a_j$ of the metric are critical points of $f$, and
$\alpha_j-1$ are the multiplicities of these critical points.

The following results are known for this case.

First of all, the sum of $\alpha_j-1$
must be even: if $d$ is the degree of $f$, then
\begin{equation}\label{1}
2+\sum_{j=0}^{n-1}(\alpha_j-1)=2d.
\end{equation}
This is stronger than the
general necessary condition (\ref{curv}).

Second,
\begin{equation}\label{2}
\alpha_j\leq d\quad\mbox{for all}\quad j,
\end{equation}
because a rational function of
degree $d$ cannot have a point where the local  degree is greater than $d$.

Subject to these two restrictions (\ref{1}) and (\ref{2}), a rational
function with prescribed critical points $a_j$
of multiplicities $\alpha_j-1$
always exists \cite{G,S}. Thus
{\em for any $a_0,\dots,a_{n-1}$ and any $\alpha_0,\dots,\alpha_{n-1}$ satisfying (\ref{1}) and (\ref{2})
there exist metrics of curvature $1$ on $S$ with angles
$2\pi\alpha_j$ at $a_j$.}

Two rational functions $f_1$ and $f_2$
are called {\em equivalent}
if $f_1=\phi\circ f_2$, where $\phi$ is a fractional-linear
transformation. Equivalent functions have the same critical points
with the same multiplicities. Equivalent functions generate equivalent metrics.

The number of equivalence classes of rational functions with prescribed
critical points and multiplicities is at most
$\K(\alpha_0-1,\ldots,\alpha_{n-1}-1)$,
where $\K$ is the {\em Kostka number} which can be described as follows.
Consider Young diagrams of shape $2\times(d-1)$ consisting of two rows
of length $d-1$. A semi-standard Young tableau (SSYT)
is a filling of such a diagram with positive integers such that an integer $k$
appears $\alpha_{k-1}-1$ times,
the entries are strictly increasing in the columns
and non-decreasing in the rows. The number of such SSYT's is the
Kostka number $\K(\alpha_0-1,\ldots,\alpha_{n-1}-1)$.

For a generic choice of the critical points $a_j$, the number
of classes of rational functions is equal to the Kostka number, see
\cite{S,EGSV}.

Suppose now that the points $a_j$ and the corresponding multiplicities
$\alpha_j$
are symmetric with respect to some circle.
We may assume without loss of generality that this circle
is the real line $\R\cup\{\infty\}$.
It may happen that among the rational functions $f$
with these given
critical points and multiplicities none is symmetric. So the resulting
metrics are all asymmetric as well \cite{EG,EG2}.

However, there is a surprising result \cite{EG0,EG,EGSV,MTV}
which was conjectured
by B. and M. Shapiro:
\vspace{.1in}

{\em If all critical points of a rational function lie on a circle,
then the function is equivalent to a function symmetric with respect to this circle. Moreover,
in this case there
are exactly $\K(\alpha_0-1,\ldots,\alpha_{n-1}-1)$ classes of rational functions
with prescribed critical points of the multiplicities $\alpha_0-1,\ldots,\alpha_{n-1}-1$.}
\vspace{.1in}

It is interesting to find out which of these results
can be extended to the general case of arbitrary positive $\alpha_j$.

Real rational functions with prescribed real critical points are classified
by combinatorial objects which are called {\em nets}.
In sections \ref{intro-to-part-2}-\ref{sec:classification}
we will use similar nets to classify spherical quadrilaterals.

\section{The cases $n\leq 3$}\label{cases-leq-3}

We recall some known results for metrics with at most $3$ conic singularities.
First of all, there is no metric on the sphere with one conic singularity.
The case of two singularities is easy: the necessary and sufficient condition
for the existence of the metric is that the angles at the singularities are
equal, and the metric with prescribed angles is unique \cite{Troy2}.

For the case of three singularities, a complete description was
obtained in \cite{UY,E,FKKRUY}:

{\em If none of the $\alpha_j$ is an integer, then the necessary
and sufficient condition for the existence of the metric is
\begin{equation}\label{3non}
\cos^2\pi\alpha_0+\cos^2\pi\alpha_1+\cos^2\pi\alpha_2+
2\cos\pi\alpha_0\cos\pi\alpha_1\cos\pi\alpha_2<1,
\end{equation}
which is equivalent to
\begin{eqnarray*}
&&\cos\pi\frac{\alpha_0+\alpha_1+\alpha_2}{2}\\
&\times& \cos\pi\frac{-\alpha_0+\alpha_1+\alpha_2}{2}
\cos\pi\frac{\alpha_0-\alpha_1+\alpha_2}{2}
\cos\pi\frac{\alpha_0+\alpha_1-\alpha_2}{2}<0,
\end{eqnarray*}
and the metric with prescribed angles is unique.

If $\alpha_0$ is an integer, then the necessary and sufficient condition
for the existence of a metric is that either $\alpha_1+\alpha_2$
or $|\alpha_1-\alpha_2|$ is an integer $m$ of the opposite parity to $\alpha_0$,
and $m\leq\alpha_0-1$. All metrics with given angles are equivalent.}
\vspace{.1in}

\begin{rmk}\label{rotation} {\rm Notice that for the case when all three $\alpha_j$
are non-integers, the condition of existence of a spherical triangle
with angles $\pi\alpha_j$ coincides with the condition on rotation angles
of three elliptic transformations to be simultaneously conjugate to rotations (see \cite{B}).}
\end{rmk}

\section{The case $n=4$ with two integer corners: condition on the angles}\label{cond-angles}

In the rest of the paper, 
we study the case $n=4$ with two integer
$\alpha_j$. We answer the following questions:

a) In the equation (\ref{heun}), for which $\alpha_j$ one can choose
$\lambda$ so that the monodromy group is conjugate to a subgroup of $PSU(2)$?

b) If $\alpha_j$ satisfy a), how many choices of $\lambda$ are
possible?

c) If, in addition, all $a_j$ are real, how many choices of real $\lambda$
are possible?

One cannot have exactly one non-integer $\alpha_j$. Indeed, in this case
the developing map $f$ will have just one branching point on the sphere,
which is impossible by the Monodromy Theorem.

Let us consider the case of two non-integer $\alpha_j$. In this section we obtain
a necessary and sufficient condition on the angles for this case, that is, solve the
problem a).

We place the two singularities corresponding to non-integer $\alpha$ at
$a_0=0$ and $a_3=\infty$,
and let the total angles at these points be $2\pi\alpha_0$ and
$2\pi\alpha_3$,
where $\alpha_0$ and $\alpha_3$ are not
integers. Then the developing map has an analytic continuation
in $\C^*$ from which we conclude that the monodromy group must
be a cyclic group
generated by a rotation $z\mapsto ze^{2\pi i\alpha}$,
with some $\alpha\in(0,1)$. This means that $f(z)$ is multiplied by
$e^{2\pi i\alpha}$ when $z$ describes a simple loop around the origin.
Thus $g(z)=z^{-\alpha}f(z)$ is a single valued function with at most
power growth at $0$ and $\infty$.
Then we have a representation
$f(z)=z^\alpha g(z)$, where $g$ is a rational function.
Then $\alpha_0=|k+\alpha|,\; \alpha_3=|j+\alpha|,$
where $k$ and $j$ are integers, so
either $\alpha_0-\alpha_3$ or
$\alpha_0+\alpha_3$ is an integer.
The angles $2\pi\alpha_1$ and $2\pi\alpha_2$ at the other two singular points
$a_1$ and $a_2$ of the metric
are integer multiples of $2\pi$, and they are the critical points of $f$
other than $0$ and $\infty$.

Let $g=P/Q$ where $P$ and $Q$ are polynomials without common zeros
of degrees $p$ and $q$, respectively.
Let $p_0$ and $q_0$ be the multiplicities
of zeros of $P$ and $Q$ at $0$. Then $\min\{ p_0,q_0\}=0$, because the fraction $P/Q$ is irreducible.

The equation for the critical points of $f$ is the following:
\begin{equation}\label{crit}
z(P'(z)Q(z)-P(z)Q'(z))+\alpha P(z)Q(z)=0.
\end{equation}
Since $\alpha_1$ and $\alpha_2$ are integers,
we have the following system of equations:
\begin{eqnarray}\label{system}
\alpha_0&=&|p_0-q_0+\alpha|,\nonumber\\
\alpha_1+\alpha_2-2&=& p+q-\max\{ p_0,q_0\},\\
\alpha_3&=&|p-q+\alpha|.\nonumber
\end{eqnarray}
The first and the last equations follow immediately from the representation
$f(z)=z^{\alpha}P(z)/Q(z)$ of the developing map.
The second equation holds because the
left-hand side of (\ref{crit}) is a polynomial of degree exactly $p+q$,
therefore
the sum of the multiplicities of its zeros $a_1$ and $a_2$
must be $p+q-\max\{ p_0,q_0\}$

Solving this system of equations (\ref{system}) in non-negative integers
satisfying $\min\{ p_0,q_0\}=0$, $p_0\leq p$, $q_0\leq q$,
we obtain the necessary and sufficient conditions
the angles should satisfy, which we state as

\begin{thm}\label{theorem1}Suppose that four points $a_0,\ldots,a_3$ on the Riemann
sphere and numbers $\alpha_j>0$, $0\leq j\leq 3$, are such that
$\alpha_1$ and $\alpha_2$ are integers $\geq 2$.

The necessary and sufficient conditions
for the existence of a metric of curvature $1$ on the sphere, with conic
singularities at $a_j$ and angles $2\pi\alpha_j$ are the following:

a) If $\alpha_1+\alpha_2+[\alpha_0]+[\alpha_3]$ is even,
then $\alpha_0-\alpha_3$ is an integer, and
\begin{equation}\label{ineq1}
\left|\alpha_0-\alpha_3\right|+2\leq \alpha_1+\alpha_2.
\end{equation}

b) If $\alpha_1+\alpha_2+[\alpha_0]+[\alpha_3]$ is odd,
then $\alpha_0+\alpha_3$ is an integer, and
\begin{equation}\label{ineq2}
\alpha_0+\alpha_3+2\leq\alpha_1+\alpha_2.
\end{equation}
\end{thm}

{\em Sketch of the proof.} For a complete proof see \cite{EGT}.
Conditions a) and b)
are necessary and sufficient for the existence of a unique solution
$p,q,p_0,q_0,\alpha$
of the system (\ref{system}) satisfying
$$\min\{ p_0,q_0\}=0,\quad p_0\leq p,\quad q_0\leq q,\quad
\alpha\in(0,1).$$
Thus the necessity of these conditions follows from our arguments above.

We may assume without loss of generality that $a_0=0,\; a_3=\infty,\; a_1=1$
and $a_2=a\in\C$.

Then we set $R(z)=z^{\max\{p_0,q_0\}}(z-1)^{\alpha_1}(z-a)^{\alpha_2}.$
The second equation in (\ref{system}) gives $\deg R=p+q$.
Now we consider the equation
\begin{equation}\label{crit2}
z(P'Q-PQ')+\alpha PQ=R.
\end{equation}
This equation must be solved in polynomials $P$ and $Q$ of degrees $p$ and $q$
having zeros of multiplicities $p_0$ and $q_0$ at $0$.
Non-zero polynomials of degree at most $p$ modulo proportionality
can be identified with the points of the complex projective space $\P^p$.
The map
$$W_\alpha:\P^p\times\P^q\to\P^{p+q},\quad (P,Q)\mapsto z(P'Q-PQ')+\alpha PQ$$
is well defined. It is a finite map between compact algebraic varieties,
and it can be represented as a linear projection of the Veronese variety.
Its degree is equal to the degree
\begin{equation}
{p+q}\choose{p}
\end{equation}
of the Veronese variety.
Thus the equation (\ref{crit2}) always has a complex solution $(P,Q)$.
The function $f=z^\alpha P/Q$ is then a developing map with the
required properties.
So conditions a) and b) are sufficient.
This completes the proof of the first statement.


It will be convenient to introduce new parameters instead of $\alpha_j$.
Besides other advantages, we eliminate the additional parameter $A$, (see
(\ref{heun}), (\ref{cond})), and
the new parameters allow us to treat the cases
a) and b) in Theorem~\ref{theorem1}
simultaneously. We will rewrite (\ref{heun}) as
\begin{equation}\label{newheun}
z(z-1)(z-a)\left(y^{\prime\prime}-\left(\frac{\sigma}{z}+\frac{m}{z-1}+
\frac{n}{z-a}\right)y^\prime\right)+\kappa(\sigma+1+m+n-\kappa)zy=\lambda y,
\end{equation}
where $\kappa$ is an integer, $\sigma\in\R$ is not an integer,
\begin{equation}\label{cc}
[\sigma]\geq -1\quad\mbox{in Case a),\; and }\quad [\sigma]<-1\quad\mbox{in Case b)}.
\end{equation}
To achieve this we put $m=\min\{\alpha_1,\alpha_2\}-1,$ and
$n=\max\{\alpha_1,\alpha_2\}-1$.

In case a), we set
$$\sigma=\min\{\alpha_0,\alpha_3\}-1
$$
and define $\kappa$ by
\begin{equation}
\label{k1}
2\kappa=-|\alpha_0-\alpha_3|+\alpha_1+\alpha_2-2.
\end{equation}
In case b), we set
$$
\sigma=-\min\{\alpha_0,\alpha_3\}-1
$$
and define $\kappa$ by
\begin{equation}
\label{k2}
2\kappa=-\alpha_0-\alpha_3+\alpha_1+\alpha_2-2.
\end{equation}
In both cases $\kappa$ is an integer because the sums in the right-hand sides
of (\ref{k1}) and (\ref{k2}) are even. Inequalities (\ref{ineq1})
and (\ref{ineq2}) give that $\kappa\geq 0$ in both cases.
Since $m+n=\alpha_1+\alpha_2-2$, we also get
\begin{equation}
\label{kmn}
2\kappa\leq m+n.
\end{equation}
Furthermore in case b),
\begin{equation}
\label{ineqc}
\sigma+1\geq\kappa-(m+n)/2
\end{equation}
because $\alpha_0+\alpha_3\geq 2\min\{\alpha_0,\alpha_3\}$. Notice that
in case a), inequality (\ref{ineqc}) holds trivially because $\sigma+1>0$.

To summarize, the new parameters are three integers $m,n,\kappa$, and one
real non-integer number $\sigma$, subject to the conditions (\ref{ineqc}) and
\begin{equation}
\label{conditions1}
0\leq m\leq n,\quad 0\leq 2\kappa\leq m+n,
\end{equation}
Parameters $\alpha_j$ are recovered by the formulas
$$(\alpha_0,\alpha_3)=(|\sigma+1|,|m+n+\sigma+1-2\kappa|),$$
$$(\alpha_1,\alpha_2)=(m+1,n+1),$$
up to a permutation of $\alpha_1$ and $\alpha_2$, and a permutation
of $\alpha_0$ and $\alpha_3$,
and Heun's equation is as (\ref{newheun}).

\section{Counting solutions}\label{counting-sol}

In most cases, there is no uniqueness in Theorem~\ref{theorem1}.
In this section we determine the number of equivalence
classes of metrics for given $a_j$ and $\alpha_j$, assuming that
two of the $\alpha_j$ are integers. As explained in the previous section,
in this case the Heun equation has a polynomial solution, and
a solution of the form $z^\alpha P$, where $P$ is a polynomial.
We call functions of this last type {\em quasipolynomials}.

Substituting a formal power series
\def\Z{{\mathbb{Z}}}
$$H(z)=\sum_{s\in\Z+\beta}h(s)z^s$$
to the equation (\ref{newheun}), we obtain recurrence relations of the form
\begin{equation}\label{rec}
c_{s-1}h(s-1)+a_s h(s)+b_s h(s+1)=0,
\end{equation}
which can be visualized as a multiplication of the vector $(h(s))$ by a
Jacobi (three-diagonal) matrix
\begin{equation}\label{matrix}
\left(\begin{array}{ccccccc}
&&\ldots&&\ldots&&\\
\ldots&c_{s-1}&a_s&b_s&0&0&\ldots\\
\ldots&0&c_s&a_{s+1}&b_{s+1}&0&\ldots\\
\ldots&0&0&c_{s+1}&a_{s+2}&b_{s+2}&\ldots\\
&&\ldots&&\ldots&&
\end{array}\right).
\end{equation}
The explicit expressions are
\begin{eqnarray}\label{coefficients1}
b_{s}&=&a(s+1)(s-\sigma),\\ \label{coefficients2}
a_{s}&=&-s\left((a+1)(s-1-\sigma)-ma-n\right)-\lambda,\\
c_s&=&(s-\kappa)(s+\kappa-\sigma-m-n-1).\label{coefficients3}
\end{eqnarray}
We see that an ``eigenvector'' $h(s)$ can be a finitely supported
sequence, say with support $[s_1,s_2]$, if and only if $b_{s_1-1}=0$
and $c_{s_2}=0$. If $b_{s_1-1}=0$, and $b_s\neq0$ for $s> s_1$,
the elements $h(s)$ can be defined recursively from (\ref{rec}),
with arbitrary non-zero value of $h(s_1)$.
Each $h(s)$ is a polynomial in $\lambda$ of degree $s-s_1$, and
the condition of termination of the sequence at the place $s_2$ is
\begin{equation}\label{co}
c_{s_2-1}h(s_2-1)+a_{s_2} h(s_2)=0.
\end{equation}
This is a polynomial equation of degree $s_2-s_1+1$ in $\lambda$ which is
the condition of having a polynomial solution of degree $d=s_2-s_1$.
Similar condition gives the existence of a quasipolynomial solution.
Solutions of (\ref{co}) are eigenvalues of $(d+1)\times(d+1)$
Jacobi matrix obtained by truncating the infinite matrix (\ref{matrix})
by leaving rows and columns with indexes from $s_1$ to $s_2$.

Substituting a formal power series
$$H(z)=\sum_{s=0}^\infty h(s) (z-1)^s$$
to the equation (\ref{newheun}) we
obtain another recurrence relation of the form (\ref{rec}) with appropriate coefficients $a_s,b_s$ and $c_s$. Since the exponents of (\ref{newheun}) at the
point $1$ are $0$ and $m+1$, we can always find a holomorphic function $H$
whose power series begins with the term $(z-1)^{m+1}$. But to find a power
series solution beginning with a constant term, the condition
\begin{equation}\label{coco}
c_{m-1}h(m-1)+a_{m}h(m)=0
\end{equation}
must be satisfied, and this is a polynomial equation of degree $d=m+1$ in
$\lambda$. If the condition (\ref{coco}) is satisfied, then $h(m+1)$
can be chosen arbitrarily. Equation (\ref{coco}) is the triviality condition
of the monodromy around $1$. Solutions of (\ref{coco}) are
eigenvalues of the $d\times d$ Jacobi matrix obtained by taking rows and
columns of (\ref{matrix}) with indices from $0$ to $m$.

Thus we have four polynomial conditions which are necessary
for Heun's equation (\ref{newheun})
to have two solutions: a polynomial and a quasipolynomial.
\vspace{.1in}

\noindent
(i) $C_1(\lambda)=0$ iff there exists a polynomial solution,
\vspace{.1in}

\noindent
(ii) $C_2(\lambda)=0$ iff there exists a quasipolynomial solution,
\vspace{.1in}

\noindent
(iii) $C_3(\lambda)=0$ iff the monodromy at $1$ is trivial, and
\vspace{.1in}

\noindent
(iv) $C_4(\lambda)=0$ iff the monodromy at $a$ is trivial.
\vspace{.1in}

The degrees of these
equations are:
$$\deg C_1=\kappa+1,\quad\deg C_2=m+n-\kappa+1,\quad
\deg C_3=m+1,\quad\deg C_4=n+1.$$
The first two formulas follow by setting $c_{s_2}=0$ in (\ref{coefficients3}),
and for the other two one has to rewrite (\ref{newheun}) to place a singular point
with integer exponents at $0$, and write the formula for $b(s)$ for this transformed
equation (see (\ref{111}) and (\ref{xxx}) below).

Thus the number of values of $\lambda$ for which all four polynomials
$C_i$ vanish is at most $\min\{\kappa+1,m+1,n+1\}$, where we used (\ref{kmn}).
Expressing this in terms of the original exponents $\alpha_j$ with the help
of (\ref{k1}) and (\ref{k2}) we obtain

\begin{thm}\label{theorem2} The number of classes of
metrics with prescribed
angles $2\pi\alpha_j$ at the given points $a_j$ is at most
$$\min\{\alpha_1,\alpha_2, \kappa+1\},$$
where $\kappa$ is defined by (\ref{k1}), (\ref{k2}),
$$
\kappa+1=\left\{\begin{array}{ll}
(\alpha_1+\alpha_2-|\alpha_0-\alpha_3|)/2& \mbox{in case a)},\\
&\\
(\alpha_1+\alpha_2-\alpha_0-\alpha_3)/2& \mbox{in case b)}.
\end{array}\right.$$
\end{thm}

We will later see that equality holds for generic $a$.
The crucial fact is

\begin{prop}\label{proposition1} Of the four polynomials
$C_j, 1\leq j\leq 4$, the polynomial of the smallest degree divides
each of the other three polynomials.
\end{prop}
\vspace{.1in}

{\em Proof}. We write Heun's equation in the form (\ref{newheun}).

1. Suppose that $C_1$ is the polynomial of the smallest degree $\kappa+1$.
For every root $\lambda$ of $C_1$,
we have a polynomial solution $p$ of degree at most $\kappa$.
So $p$ cannot have a zero of
order $\alpha_1$ or $\alpha_2$, because by assumption
these numbers are at least $\kappa+1$. This implies
that the monodromy is trivial at $1$ and $a$. Then the second solution
of Heun's equation also has no singularities at $1$ and $a$.
So in this case $C_1$ divides $C_2$, $C_3$ and $C_4$.

2. $C_2$ cannot be the polynomial of the smallest degree in view of
(\ref{kmn}).

3. It remains to show that if $C_3$ is the polynomial of the smallest degree
and $\lambda$ is a root of $C_3$, then $\lambda$ is also a root of $C_1$
and $C_2$. Then it will follow that there is a polynomial and a quasipolynomial
solutions, so the monodromy will be also trivial at $a$, that is, $\lambda$
will be a root of $C_1$, $C_2$ and $C_4$.
The situation in the case when $C_4$
is of the smallest degree is completely similar. The result follows from

\begin{prop}\label{proposition2}
 Let $m,n,\kappa$ be integers,
$0\leq m\leq n,\; 0\leq 2\kappa\leq m+n$, and $\sigma$ not an integer.
Consider the differential equation (\ref{newheun})
which we write as $Dy=0$, where
$$Dy = z(z-1)(z-a)\left(y''-\left(\frac{\sigma}{z}+\frac{m}{z-1}+\frac{n}{z-a}
\right)y'\right)+\kappa(\sigma+1+\tau)zy-\lambda y,$$
and $\tau=m+n-\kappa\geq \kappa$. Suppose that $m\leq \kappa$, and that the monodromy at $1$
is trivial, that is
all solutions are holomorphic at the point $1$. Then there exist
a polynomial solution and a quasipolynomial solution.
\end{prop}
\vspace{.1in}

{\em Proof.}
Let us transform our equation (\ref{newheun}) to the form
\begin{eqnarray}\displaystyle\nonumber
&&z(z-1)(z-a)\left(y^{\prime\prime}-
\left(\frac{m}{z}+\frac{\sigma}{z-1}+
\frac{\sigma+2+m+n-2\kappa}{z-a}\right)y'\right)\\
&&+
\kappa(\kappa-n-1)zy=\lambda y.\label{111}
\end{eqnarray}
It has a polynomial solution of degree $k$ simultaneously with the
original equation (\ref{newheun}). Consider the infinite Jacobi matrix
(\ref{matrix}). The triviality of the monodromy of (\ref{111}) at $0$
means that $\lambda$ is an eigenvalue of the truncated Jacobi matrix $J_0$ given by the first $m+1$ rows and columns. The existence of a polynomial solution
of (\ref{111}) of degree $k$ means that $\lambda$ is an eigenvalue of the
truncated matrix $J_1$ given by the first $k+1$ rows and columns.

By explicit formulas for the entries, the matrix $J_1$ is upper
block-triangular, the bottom-left $(m+1)\times k$ block being equal to zero,
and the top-left $(m+1)\times (m+1)$ block of $J_1$ equals $J_0$.
Thus every
eigenvalue of $J_0$ is an eigenvalue of $J_1$.

To show the existence of a quasipolynomial solution $y(z)=z^{\sigma+1}q(z)$,
we
write the differential equation for $q$ ($\sigma$ will be replaced
by $-\sigma-2$) and then transform it to the form (\ref{111}).
\vspace{.1in}

To summarize the contents of this section, we consider, for any
given $\alpha_0,\ldots,\alpha_3$ satisfying conditions a) or b) of
Theorem~\ref{theorem1}, the polynomial $F(a,\lambda)$ which is the
polynomial of the smallest degree of those $C_j$ in
Proposition~\ref{proposition1}. The condition
\begin{equation}\label{maineq}
F(a,\lambda)=0
\end{equation}
is equivalent to the statement that the monodromy of Heun's equation
is conjugate to a subgroup of $PSU(2)$.
Thus equivalence classes of metrics of positive curvature $1$ with
singularities at $0,1,a,\infty$ with prescribed $\alpha_j$ are in one-to-one
correspondence with solutions of the equation (\ref{maineq}).

{\em Remark.} The value $\lambda$ in (\ref{maineq}) depends not only
on the quadrilateral (or a metric) that we consider but also on
the choice of the Heun equation. Different Heun equations corresponding
to the same quadrilateral can be obtained by cyclic permutation of the 
vertices, and by the different choices of exponents at the singularities.
The angle at a vertex only fixes the absolute value of the difference of the
exponents. The values of $\lambda$ corresponding to the same
quadrilateral but different Heun equations are related by
fractional-linear transformations.

\section{Counting real solutions}\label{counting-real}

In this section we assume that $a$ is real and estimate from below
the number of real Heun's equations with given conic singularities at
$0,1,a,\infty$ with prescribed angles and unitary monodromy, or,
which is the same, the number of real solutions $\lambda$ of equation
(\ref{maineq}).
We will also show that for generic $a$ we have equality
in the inequality of Theorem~\ref{theorem2}
for the number of complex solutions.

Our estimates will be based on the following lemma.

\begin{lemma}\label{lemma1} Let $J$ be a real $(d+1)\times (d+1)$ Jacobi matrix
$$J=\left(\begin{array}{cccccccc}
a_0&b_0&0&0&\ldots&0&0&0\\
c_0&a_1&b_1&0&\ldots&0&0&0\\
0&c_1&a_2&b_2&\ldots&0&0&0\\
\ldots&\ldots&\ldots&\ldots&\ldots&\ldots&\ldots&\ldots\\
\ldots&\ldots&\ldots&\ldots&\ldots&\ldots&\ldots&\ldots\\
0&0&0&0&\ldots&c_{d-2}&a_{d-1}&b_{d-1}\\
0&0&0&0&\ldots&0&c_{d-1}&a_d
\end{array}\right).
$$
a) If $b_jc_j>0$, $0\leq j\leq d-1$, then all eigenvalues of $J$ are real
and simple.
\vspace{.1in}

\noindent
b) If $c_j\neq 0$ for $0\leq j\leq d-1$, then we have
\begin{equation}\label{herm}
J^TR=RJ,
\end{equation}
\def\diag{{\mathrm{diag}}}
where $R=\diag(r_0,\ldots,r_d),$ $r_0=1$ and
$$r_j=r_{j-1}\frac{b_{j-1}}{c_{j-1}},\quad 1\leq j\leq d.$$

\noindent
c) Suppose that the sequence $d_j=b_jc_j$ has the property
$d_j>0$ for $0\leq j\leq k$ and $d_j<0$ for $k+1\leq j\leq d$.
Then the number of pairs of
non-real eigenvalues, counting multiplicity,
is at most $[(d-k)/2]$.
\end{lemma}

{\em Proof.} Statement a) in contained in \cite{GK}; we include a simple
proof for convenience.
\def\diag{\mathrm{diag}}
Consider the matrix $S=\diag(s_0,\ldots,s_d),\; s_0=1$, and
$$s_j=s_{j-1}\sqrt{b_{j-1}/c_{j-1}},\quad 1\leq j\leq d.$$
Under the assumption of part a), the fraction under the square root
is positive, and the matrix $S$ is the positive square root of the
matrix $R$ given in part b), $S^2=R$. By explicit calculation, the matrix
$$\tilde{J}=SJS^{-1}$$
is real and symmetric, so it is diagonalizable and has real eigenvalues.
The top right $d\times d$ submatrix of $\tilde{J}-\lambda I$,
where $I$ is the identity matrix and $\lambda$ is an eigenvalue of $\tilde{J}$,
is lower triangular and has the determinant
$$\sqrt{b_0\ldots b_{d-1}c_0\ldots c_{d-1}}\neq 0.$$
Hence all eigenvalues of $\tilde{J}$ are simple.

Statement b) is proved by direct calculation.

To prove statement c), we notice that our assumption about signs of $d_j$
implies that there are $[(d-k)/2]$ negative numbers among
$r_1,\ldots,r_n$, and the rest are positive.
Condition (\ref{herm}) means that our matrix $J$
is symmetric with respect to the bilinear form $(x,y)_R=x^TRy$,
that is
$$(Jx,y)_R=(x,Jy)_R.$$
Quadratic form $(x,x)_R$ has $[(d-k)/2]$ negative squares,
so our matrix $J$ has at most $2[(d-k)/2]$ non-real eigenvalues,
counted with algebraic multiplicities, according to the theorem
of Pontrjagin \cite{Pnt}. This proves the lemma.

Our main result on the real case is the following:

\begin{thm}\label{theorem3} Consider the metrics of curvature $1$
on the sphere with real conic
singularities $a_0,a_1,a_2,a_3$ and the corresponding angles $2\pi\alpha_0,
2\pi\alpha_1,2\pi\alpha_2,$ $2\pi\alpha_3$, where $\alpha_1$ and $\alpha_2$
are integers. Suppose that conditions of Theorem~\ref{theorem1} are satisfied. Then
\vspace{.1in}

\noindent
(i) If the pairs $(a_0,a_3)$ and $(a_1,a_2)$ do not separate each
\def\bR{{\mathbf{\overline{R}}}}
other on the circle $\bR$, then all metrics with these
angles and singularities are symmetric. Their number
is equal to
$$\min\{\alpha_1,\alpha_2,\kappa+1\}$$
where $\kappa$ is defined in (\ref{k1}) and (\ref{k2}).
\vspace{.1in}

\noindent
(ii) If the pairs $(a_0,a_3)$ and $(a_1,a_2)$ separate each other, then
the number of of classes of symmetric metrics is at least
$$\min\{\alpha_1,\alpha_2,\kappa+1\}-2\left[\frac{1}{2}
\min\left\{\alpha_1,\alpha_2,\delta\right\}\right],$$
where
$$\delta=\frac{1}{2}\max\{
\alpha_1+\alpha_2-[\alpha_0]-[\alpha_3], 0\}
$$

\noindent
(iii) There is an $\epsilon>0$ depending on the $\alpha_j$ such that
if
\begin{equation}\label{crossrat}
\left|\frac{(a_2-a_0)(a_3-a_1)}{(a_1-a_0)(a_3-a_2)}\right|<\epsilon
\end{equation}
in (\ref{heun}) then all of metrics are symmetric,
and their number is as in (i).
\end{thm}

The expression under the absolute value
in the left hand side of (\ref{crossrat}) is the
cross-ratio which is equal to $a$ when the vertices
are $0,1,a,\infty$.

In section \ref{sec:examples} we will show that the estimate in (ii)
is achieved sometimes,
and in section \ref{alternative} we give another independent proof
of this estimate.

\begin{cor}\label{corollary1} Suppose that $a_j$ are real,
If the pairs $(a_0,a_3)$ and $(a_1,a_2)$
separate each
other on $\bR$, and
$$[\alpha_0]+[\alpha_3]+2\geq \alpha_1+\alpha_2$$
then all metrics with singularities at $a_j$ and angles $\pi\alpha_j$
are symmetric with respect to the real line.
\end{cor}

If Case a) of Theorem \ref{theorem1} prevails, this corollary can be also obtained
as a special case of Theorem 5.2 from \cite{MTV}.
\vspace{.1in}

{\em Proof of Theorem \ref{theorem3}.}
Let us transform our equation (\ref{newheun}) to the form
(\ref{111}).
This equation has the same exponents at the singularities as (\ref{newheun}),
but we placed the point with the smaller integer exponent $m$ at $0$.
The recurrence relations similar to (\ref{rec}) have in this case
the following coefficients
\begin{eqnarray}\label{xxx}
c_{s}&=&(\kappa-s)(\kappa-n-s-1),\nonumber\\
a_{s}&=&-s(s+\sigma+1+n-2\kappa+a(s-1-m-\sigma))\\
b_{s}&=&a(s+1)(s-m).\nonumber
\end{eqnarray}
So
\begin{equation}\label{po1}
b_sc_s=(\kappa-s)(\kappa-n-s-1)a(s+1)(s-m),
\end{equation}
which is positive for $a>0$ and
$0\leq s<\min\{\kappa,m\}$. Thus for $a>0$ all eigenvalues are real and distinct
by Lemma~\ref{lemma1} a) with $d=\min\{m,\kappa\}+1$.
This proves (i).

For the case (ii) that is $a<0$, we transform
equation (\ref{111}) by the change of the variable $z'=1-z$ into equation
\begin{eqnarray}\nonumber\displaystyle
&&z(z-1)(z-a')\left( y^{\prime\prime}-
\left(\frac{\sigma}{z}+\frac{m}{z-1}+\frac{\sigma+m+n+2-2\kappa}{z-a'}\right)y'\right)\\
&&+\kappa(\kappa-n-1)zy
=\lambda y.\label{222}
\end{eqnarray}
Here $a'=1-a$. The coefficients of the recurrence become
\begin{eqnarray}
c_s&=&(\kappa-s)(\kappa-n-s-1)\nonumber\\
a_s&=&-s(s+m+n+1-2\kappa)+a'(s-m-\sigma-1)),\label{75}\\
b_s&=&a'(s+1)(s-\sigma).\nonumber
\end{eqnarray}
So we have
\begin{equation}\label{po2}
c_sb_s=(\kappa-s)(\kappa-n-s-1)a'(s+1)(s-\sigma),
\end{equation}
which is positive when
$a'>0, \; \sigma+1>\kappa,$ and $0\leq s< \kappa$. Thus under this condition, all
eigenvalues are real. The range $a<0$ is covered by $a'>0$.

If $\sigma+1\leq \kappa$, then the Jacobi matrix with entries (\ref{75})
has the property described in
Lemma~\ref{lemma1} c), where $R$ has $[\sigma]+1$ positive squares and $[(\kappa-[\sigma])/2]$
negative squares. So the number of pairs of
non-real eigenvalues, counting algebraic multiplicities,
is at most $(\kappa-[\sigma])/2.$
We recall that $\kappa+1$ is the degree of the polynomial $C_1$ in
Proposition \ref{proposition1}. So the polynomial of minimum degree among the $C_j$
has at most $(1/2)\min\{ \kappa+1,m+1,\kappa-[\sigma]\}$ pairs of non-real zeros,
and using the value of $\kappa$ from (\ref{k1}), (\ref{k2}) and
the inequality (\ref{cc}), we obtain (ii).

To prove (iii) we notice that when $a=0$ in (\ref{xxx}), the Jacobi matrix
is triangular, so its eigenvalues are real and simple. By continuity
this situation persists when $|a|$ is small enough.

\section{Introduction to nets}\label{intro-to-part-2}

In this section we begin a different treatment of spherical polygons
which is independent of sections \ref{DE-connection}--\ref{counting-real}.

\begin{df}\label{marked} {\rm A spherical $n$-gon $Q$ is {\em marked} if
one of its corners, labeled $a_0$, is identified as the first corner, and
the other corners are labeled so that
$a_0,\ldots,a_{n-1}$ are in the counterclockwise order on the boundary of $Q$.}
\end{df}

We call $Q$ a {\em spherical polygon} when $n$ is not specified.
When $n=2,\,3$ and $4$, we call $Q$ a spherical {\em digon}, {\em triangle}
and {\em quadrilateral}, respectively.
For $n=1$, there is a unique marked 1-gon with the angle $\pi$ at its single corner.
For convenience, we often drop ``spherical'' and refer simply to $n$-gons,
polygons, etc.

Let $Q$ be a marked spherical polygon and $f:Q\to\S$ its developing map.
The images of the sides $(a_j,a_{j+1})$ of $Q$ are contained in
geodesics (great circles) on $\S$. These geodesics define a {\em partition}
$\P$ of $\S$ into vertices (intersection points of the circles) edges
(arcs of circles between the vertices) and faces (components of the complement
to the circles). Some corners of $Q$ may be integer (i.e., with angles $\pi\alpha$
where $\alpha$ is an integer). Two sides of $Q$ meeting
at its integer corner are mapped by $f$ into the same circle.

The corners of $Q$ with integer (resp., non-integer) angles are called its {\em integer}
(resp., {\em non-integer}) corners.
The {\em order} of a corner is the integer part of its angle.
A {\em removable} corner is an integer corner of order 1.
A polygon $Q$ with a removable corner is isometric to a polygon with a smaller number of corners.

A polygon with all integer corners is called {\em rational}.
All sides of a rational polygon map to the same circle.

\begin{df}\label{def:net}
{\rm Preimage of $\P$ defines a cell decomposition $\Q$ of $Q$, called the {\em net} of $Q$.
The corners of $Q$ are vertices of $\Q$.
In addition, $\Q$ may have {\em side vertices}
and {\em interior vertices}.
If the circles of $\P$ are in general position, interior vertices have degree 4,
and side vertices have degree 3.
Each face $F$ of $\Q$ maps one-to-one onto a face of $\P$.
An edge $e$ of $\Q$ maps either onto an edge of $\P$ or onto a part of an edge of $\P$.
The latter possibility may happen when $e$ has an end at an integer corner of $Q$.
The adjacency relations of the cells of $\Q$ are compatible with the adjacency
relations of their images in $\S$.
The net $\Q$ is completely defined by its 1-skeleton, a connected
planar graph. When it does not lead to confusion, we use the same
notation $\Q$ for that graph.

If $C$ is a circle of $\P$, then the intersection $\Q_C$ of $\Q$ with the preimage
of $C$ is called the $C$-{\em net} of $Q$.
Note that the intersection points of $\Q_C$ with preimages of other circles of $\P$
are vertices of $\Q_C$.
A $C$-{\em arc} of $Q$ (or simply an arc when $C$ is not specified) is
a non-trivial path $\gamma$ in the 1-skeleton of $\Q_C$
that may have a corner of $Q$ only as its endpoint.
If $\gamma$ is a subset of a side of $Q$ then it is a {\em boundary arc}.
Otherwise, it is an {\em interior arc}.
The {\em order} of an arc is the number of edges of $\Q$ in it.
An arc is {\em maximal} if it is not contained in a larger arc.
Each side $L$ of $\Q$ is a maximal boundary arc.
The {\em order} of $L$ is, accordingly, the number of edges of $\Q$ in $L$.}
\end{df}

\begin{df}\label{df:primitive}
{\rm We say that $Q$ is {\em reducible} if it
contains a proper polygon with the corners at some (possibly, all) corners of $Q$.
Otherwise, $Q$ is {\em irreducible}.
The net of a reducible polygon $Q$ contains an interior arc
with the ends at two distinct corners of $Q$.
We say that $Q$ is {\em primitive} if it is irreducible and its net does not contain
an interior arc that is a loop.}
\end{df}

\begin{df}\label{def:combequiv}
{\rm Two irreducible polygons $Q$ and $Q'$ are
{\em combinatorially equivalent} if there is an orientation preserving
homeomorphism $h:Q\to Q'$ mapping the corners of $Q$ to the corners of $Q'$,
and the net $\Q$ of $Q$ to the net $\Q'$ of $Q'$.

Two rational polygons $Q$ and $Q'$ with all sides mapped
to the same circle $C$ of $\P$ are
combinatorially equivalent if there is an orientation preserving
homeomorphism $Q\to Q'$ mapping the net $\Q_C$ of $Q$
to the net $\Q'_C$ of $Q'$.

If $Q$ and $Q'$ are reducible and represented
as the union of two polygons $Q_0$ and $Q_1$ (resp., $Q'_0$ and $Q'_1$)
glued together along their common side, then $Q$ and $Q'$ are
combinatorially equivalent when there is an orientation preserving
homeomorphism $h:Q\to Q'$ inducing combinatorial equivalence between $Q_0$ and $Q'_0$,
and between $Q_1$ and $Q'_1$.

For marked polygons $Q$ and $Q'$, we require also that the marked corner $a_0$ of $Q$
is mapped by $h$ to the marked corner $a'_0$ of $Q'$.}
\end{df}

Thus an equivalence class of nets is a
combinatorial object. It is completely determined by the labeling of the corners
and the adjacency relations.
We'll call such an equivalence class ``a net'' when this would not lead to confusion.

Conversely, given labeling of the corners and a partition $\Q$
of a disk with the adjacency relations
compatible with the adjacency relations of $\P$,
a spherical polygon with the
net $\Q$ can be constructed by gluing together the cells of $\P$
according to the adjacency relations of $\Q$.
Such a polygon is unique if the image of $a_0$, the direction in which the
image of the edge $(a_0,a_1)$ is traversed, and the images of integer vertices which are different from
the vertices of $\P$, are fixed.

In what follows we classify all equivalence classes of nets in the case when $\P$ is defined
by two circles. In this case, the boundary of each 2-cell of the net $\Q$ of $Q$ consists
of two segments mapped to the arcs of distinct circles, with the vertices at the common
endpoints of the two segments and, possibly, at some integer corners of $Q$.

By the Uniformization Theorem, each marked spherical polygon is conformally equivalent
to a closed disk with marked points on the boundary.
In the case of a quadrilateral, we have four marked points, so conformal class
of a quadrilateral depends on one parameter, the modulus of the quadrilateral.
In section \ref{chains} below
we will study whether for given permitted angles of a quadrilateral an arbitrary modulus
can be achieved. This will be done by the method of continuity,
and for this we'll need some facts about deformation of spherical quadrilaterals
(see section \ref{deformation}).

\section{Nets for a two-circle partition}\label{sec:intro}
\begin{figure}
\centering
\includegraphics[width=3in]{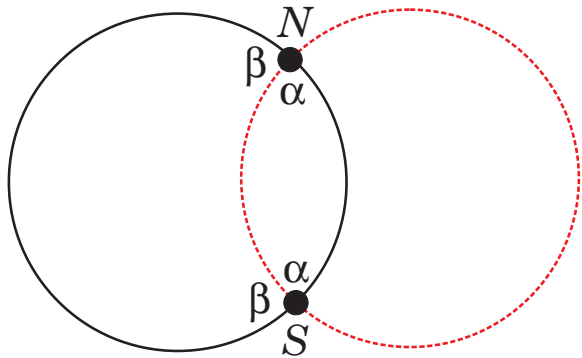}
\caption{Partition $\P$ of the Riemann sphere by two circles.}\label{partition}
\end{figure}

Let us consider a partition $\P$ of the Riemann sphere $\S$
by two transversal circles intersecting at the angle $\alpha$ (see Fig.~\ref{partition}).
{\em Vertices} $N$ and $S$ of $\P$ are the intersection points of the two circles.

We measure the angles in multiples of $\pi$, so that $0<\alpha<1$,
and the complementary to $\alpha$ angle is $\beta=1-\alpha$.
An angle that is an integer multiple of $\pi$ is called, accordingly, an {\em integer angle}.
A corner with an integer angle is called an {\em integer corner}.

Let $Q$ be a spherical $n$-gon over $\P$ (see Definition \ref{n-gon}).
We assume $Q$ to be a {\em marked} polygon (see Definition \ref{marked}).

\begin{thm}\label{non-integer} An irreducible spherical polygon $Q$ over the partition $\P$
has at most two non-integer corners.
\end{thm}

\medskip
\noindent{\bf Proof.}
We prove this statement by induction on the number $m$ of faces of the net $\Q$
of $Q$. If $m=1$ then $Q$ is isometric to a face of $\P$, thus it has exactly two
non-integer corners.

Let $m>1$. Suppose first that $\Q$ has a maximal interior arc $\gamma$ that is not a loop.
Let $p$ and $q$ be the endpoints of $\gamma$.
Since $\gamma$ is maximal, both $p$ and $q$ are at the boundary of $Q$.
Since $Q$ is irreducible, at least one of them, say $p$, is not a corner of $Q$.
Thus $\gamma$ partitions $Q$ into two polygons, $Q'$ and $Q''$, each of them having
less than $m$ faces of its net.
The induction hypothesis applied to $Q'$ and $Q''$ implies that each of them
has at most two non-integer corners.
But the corners of $Q'$ and $Q''$ at $p$ are non-integer,
while $Q$ does not have a corner at $p$.
Thus $Q$ has at most two non-integer corners.

Consider now the case when all maximal interior arcs of $\Q$ are loops.
Since the 1-skeleton of $\Q$ is connected, there exists a maximal interior arc $\gamma$
of $\Q$ with both ends at a corner $p$ of $\Q$.
We may assume that the disk $D$ bounded by $\gamma$ does not contain
another arc of $\Q$ with both ends at $p$,
otherwise we can replace $\gamma$ by a smaller loop.
Let $C$ be the circle of $\P$ such that $\gamma$ is an arc of $\Q_C$,
and let $C'$ be the other circle of $\P$.
Then $\gamma$ intersects with $\Q_{C'}$ at exactly two points.
Otherwise, either $D$ would be a face of $\Q$ with
all its boundary in $\Q_C$, or $D$ would contain a face of $\Q$ with
more than one segment of both $\Q_C$ and $\Q_{C'}$ in its boundary.
Let $\gamma'$ be the maximal interior arc of $\Q_{C'}$
intersecting $\gamma$ at those two points.
If one of those points is $p$, then $p$ must be a preimage of a vertex of $\P$.
Then $\gamma'$ is a loop having both ends at $p$,
and a single intersection point $q$ with $\gamma$ inside $Q$.
But this is impossible because the complement to the union of the disks
bounded by $\gamma$ and $\gamma'$ would contain
a face of $\Q$ whose boundary would not be a circle.

\begin{figure}
\centering
\includegraphics[width=4.0in]{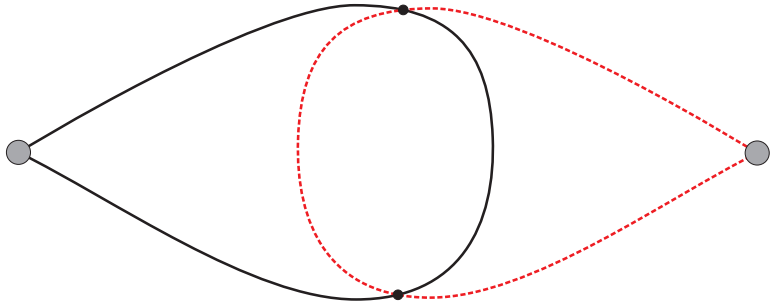}
\caption{Pseudo-diagonal connecting two integer corners of $Q$.}\label{2circles-i}
\end{figure}

Thus $p$ cannot be a preimage of a vertex of $\P$, and both intersection
points $q$ and $q'$ of $\gamma$ and $\gamma'$ are interior vertices of $\Q$.
Then $\gamma'$ must have both ends at a corner $p'$ of $\Q$.
Otherwise the complement to the union of the disks
bounded by $\gamma$ and $\gamma'$ would contain
a face of $\Q$ whose boundary would not be a circle.
The same arguments as above imply that $p'$ is not a preimage
of a vertex of $\P$, thus the union of $\gamma$ and $\gamma'$
is a {\em pseudo-diagonal} of $Q$ shown in Fig.~\ref{2circles-i}.
Removing these two loops, we obtain a polygon having
$m-4$ faces in its net, with the same number of non-integer
corners as $Q$. By the induction hypothesis, $Q$
must have at most two non-integer corners.
This completes the proof of Theorem \ref{non-integer}.

\bigskip
\noindent{\bf Rational spherical polygons.}
All corners of a rational polygon $Q$ are integer, and all its sides map to the same
circle $C$ of $\P$. Thus $Q$ is completely determined, up to combinatorial equivalence,
by its $C$-net $\Q_C$.
Each maximal arc of $\Q_C$ connects two corners of $Q$.
If $Q$ is irreducible, $\Q_C$ does not have interior arcs.
Note that converse is not true, as there may be an arc
of $\Q_{C'}$ connecting two corners of $Q$.

\section{Primitive spherical polygons with two\newline non-integer corners}\label{primitive-2}
In this section, $Q$ denotes a marked primitive spherical $n$-gon
with two non-integer corners over the partition $\P$,
one of these two corners labeled $a_0$, and the other one labeled $a_k$ where $0<k<n$.
The sides of $Q$ are labeled $L_j$ so that $a_{j-1}$ and $a_j$
are the ends of $L_j$, with $a_n$ identified with $a_0$.
We assume that the sides $L_j$ for
$1\le j\le k$ belong to the preimage of a circle $C$ of $\P$,
while the sides $L_j$ for $k<j\le n$ belong to the preimage
of the circle $C'\ne C$ of $\P$.

\begin{lemma}\label{no-interior-vertex}
The net $\Q$ of $Q$ does not have interior vertices.
\end{lemma}

\medskip
\noindent{\bf Proof.}
Let $q$ be an interior vertex of $\Q$, and let $F$ be a face of $\Q$ adjacent to $q$.
Let $\gamma$ be a maximal arc of $\Q_C$ through $q$,
where $C$ is a circle of $\P$.
Since $\gamma$ is an interior maximal arc, it may end either
at a corner of $Q$ or at a side.
Since $Q$ is primitive, $\gamma$ is not a loop, and
the ends of $\gamma$ cannot be at two distinct corners of $Q$.
If an end $p$ of $\gamma$ is at a side of $Q$, let
$L$ and $\gamma'$ be the maximal arcs of $\Q_{C'}$,
where $C'\ne C$, passing through $p$ and $q$, respectively.
Note that that $L$ is a side of $Q$ while $\gamma'$ is an interior arc.
We may assume that there are two adjacent faces, $F$ and $F'$, of $\Q$
having a common segment $pq$ of $\gamma$ in their boundary.
If there are no such faces then we can replace $q$
by an interior vertex of $\Q$ on $\gamma$ closest to $p$.
Then each of these two faces must have an integer corner
in its boundary where $L$ and $\gamma'$ intersect,
otherwise the intersection of its boundary with $\Q_C$
would not be connected.
The two corners must be distinct, since they are the ends
of a side $L$ of $Q$.
This implies that an interior arc $\gamma'$ has its ends
at two distinct corners of $Q$, thus $Q$ is not irreducible.

\begin{cor}\label{arc-order-one}
Each interior arc of the net $\Q$ of $Q$ is maximal and has order one.
\end{cor}

\begin{lemma}\label{corner-order-zero}
Each of the two non-integer corners of $Q$ has order zero.
\end{lemma}

\medskip
\noindent{\bf Proof.}
Let $p$ be a non-integer corner of $Q$ of order greater than zero.
Since $p$ is mapped to a vertex of $\P$,
there is a face $F$ of the net $\Q$ of $Q$ having $p$ as its vertex,
with two interior arcs, $\gamma$ and $\gamma'$, adjacent to $p$ in its boundary.
The arcs $\gamma$ and $\gamma'$ belong to preimages of two different circles of $\P$.
The other ends of $\gamma$ and $\gamma'$ cannot be corners of $Q$,
thus they must be side vertices of $\Q$.
This implies that the preimage of each of the two circles of $\P$ in the boundary of $F$
is not connected, a contradiction.

\begin{cor}\label{arcs}
Any interior arc of $\Q$ has one of its ends at an integer corner $a_i$
of $Q$ and another end on the side $L_j$, where either $0<i<k<j\le n$
or $0<j\le k<i<n$.
\end{cor}

\begin{df}\label{d:types} {\rm Let $Q$ be a marked primitive $n$-gon
with two non-integer corners labeled $a_0$ and $a_k$, and let $\Q$ be the net of $Q$.
For each pair $(i,j)$, let $\mu(i,j)$ be the number of interior arcs
of $\Q$ with one end at the integer corner $a_i$ and the other end on the side $L_j$.
Note that $\mu(i,j)$ may be positive only when either $0<i<k<j\le n$ or
$0<j\le k<i<n$, due to Corollary \ref{arcs}.
We call the set $\T$ of the pairs $(i,j)$ for which $\mu(i,j)>0$ the $(n,k)$-{\em type}
of $Q$ (or simply the type of $Q$ when $n$ and $k$ are fixed),
and the numbers $\mu(i,j)$ the {\em multiplicities}.
We'll see in Lemma \ref{maxtype} below that the number of pairs in
an $(n,k)$-type is at most $n-2$.
An $(n,k)$-type with exactly $n-2$ pairs is called {\em maximal}.}
\end{df}

Since interior arcs of $\Q$ do not intersect inside $Q$,
the $(n,k)$-type of $Q$ cannot contain two pairs $(i_0,j_0)$ and $(i_1,j_1)$
satisfying any of the following four conditions:
\begin{subequations}\label{forbidden}
\begin{align}
        i_0<i_1<&\;k<j_0<j_1,\\
        j_0<j_1\le&\;k<i_0<i_1,\\
        i_0<j_1\le&\;k<j_0<i_1,\\
        j_1\le i_0<&\;k<i_1<j_0.
\end{align}
\end{subequations}

Interior arcs of $\Q$ can be canonically ordered,
starting from the arc closest to the marked corner $a_0$,
so that any two consecutive arcs belong to the boundary of a cell of $\Q$.
The linear order on the interior arcs of $\Q$ induces interior
order on the pairs $(i,j)$ in the type $\T$ of $Q$.

\begin{lemma}\label{maxtype}
The number of pairs in the $(n,k)$ type $\T$ of $Q$ is at most $n-2$.
\end{lemma}

\medskip
\noindent{\bf Proof.}
Let ${i,j}\in\T$ be the first pair, corresponding to the interior arcs
of $\Q$ closest to $a_0$.
We may assume that $1\le i<k$ and $n-k<j\le n$.
Otherwise, we exchange $k$ and $n-k$.
Let $(i,m)\in\T$ be the pair furthest from $a_0$, with the same $i$ as the first pair.
Then there are at most $n-m+1$ pairs in $\T$ with the same index $i$.
The last arc of $\Q$ with the ends at the vertex $a_i$ and the side $L_m$
partitions $Q$ into two polygons, $Q'$ and $Q''$, with $Q''$ containing $a_0$.
Contracting $Q''$ to a point, we obtain a $(m-i)$-gon $\tilde Q$ with the $(m-i,k-i)$-type
$\tilde\T$ obtained from $\T$ by deleting all pairs with the same $i$ as the first one,
and relabeling vertices and sides.
Since $m\le n$ and $i>0$, we may assume inductively that the type $\tilde\T$ of $\tilde Q$
has at most $m-i-2$ pairs. This implies that the type $\T$ has at most
$n-i-1\le n-2$ pairs.

\begin{lemma}\label{non-maximal}
Any $(n,k)$-type can be obtained from a (non-unique) maximal $(n,k)$-type
if some of the multiplicities are permitted to be zero.
\end{lemma}

\medskip
\noindent{\bf Proof.}
Let $\T$ be a $(n,k)$-type with less than $n-2$ pairs.
We want to show that one can add a pair to $\T$.
We prove it by induction on $n$, the case $n=2$ being trivial.
We use notations of the proof of Lemma \ref{maxtype}.
Note first that, if $i>1$ then a pair $(1,n)$ can be added to $\T$.
Thus we may assume that $i=1$.
Next, $\T$ should contain all $n-m+1$ pairs $(1,m),\dots,(1,n)$,
otherwise a missing pair can be added to $\T$.
Finally, we can assume inductively that $\T'$ contains exactly $m-i-2=m-3$ pairs.
Thus the number of pairs in $\T$ should be $(n-m+1)+(m-3)=n-2$.

\medskip
We can associate to a maximal $(n,k)$-type $\T$ a sequence of positive integers
${\bf m}=\{m_1,\,m_2,\dots\}$ partitioned into two subsets $I$
and $J$ (we write $\underline m$ for $m\in I$ and $\overline m$ for $m\in J$) so that
\begin{equation}\label{IJ}
\sum_\nu m_\nu=n-2,\quad |I|<k,\quad,|J|<n-k.
\end{equation}
Here $m_\nu$ are the numbers of pairs $(i,j)$ in $\T$ with the same $i$,
ordered according to the linear order on the pairs in $\T$.
Obviously, for given $(n,k)$, a sequence $\bf m$ with a partition $(I,J)$ satisfying (\ref{IJ})
corresponds to at most one maximal $(n,k)$-type $\T$.

\begin{thm}\label{primitive}
For any $n\ge 2$, any $k,\;0<k<n$, any set $\T$ of pairs $(i,j)$ with either $0<i<k<j\le n$
or $0<j\le k<i<n$, such that no two pairs $(i_0,j_0)$ and $(i_1,j_1)$ in $\T$
satisfy any of the conditions (\ref{forbidden}a-d),
and any positive integers $\nu(i,j),\;(i,j)\in\T$,
there exists a marked primitive spherical $n$-gon, unique up to combinatorial
equivalence, with two non-integer corners, one of them marked, with the type $\T$
and multiplicities $\nu(i,j)$.
\end{thm}

\begin{cor}\label{ex:digon}
Each primitive digon ($n$-gon for $n=2$) with two non-integer corners
maps one-to-one to a face of $\P$,
with its two corners mapped to distinct vertices of $\P$.
There is a single (empty) $(2,1)$-type of a primitive digon.
\end{cor}

\begin{figure}
\centering
\includegraphics[width=3.0in]{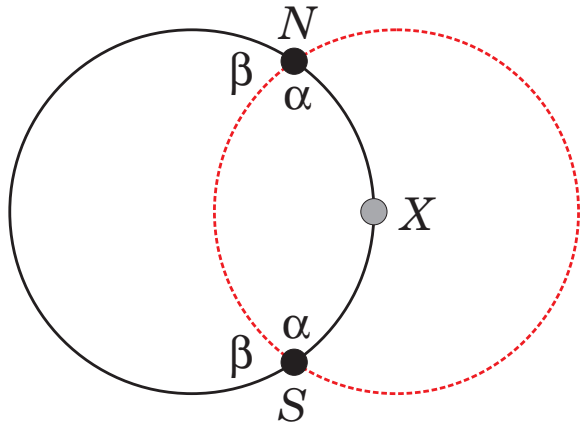}
\caption{Location of the point $X$.}\label{2circles-d}
\end{figure}

\begin{figure}
\centering
\includegraphics[width=4.0in]{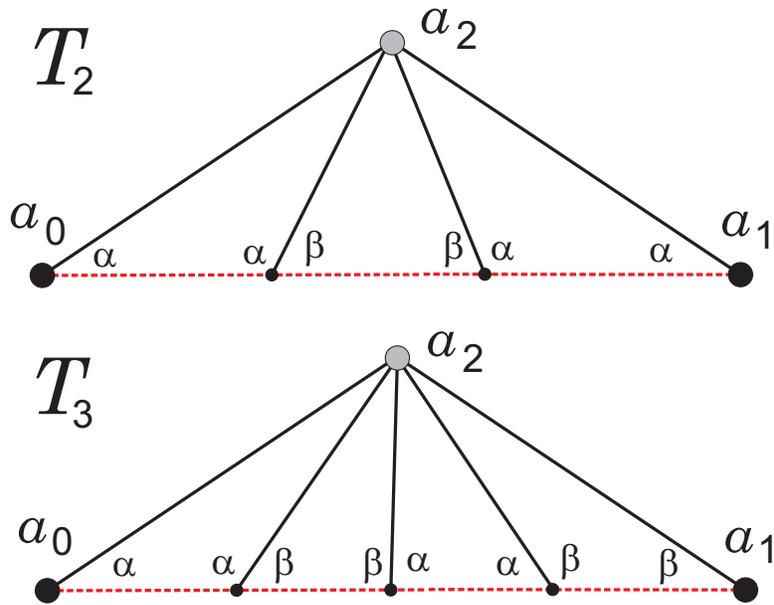}
\caption{Primitive triangles $T_\mu$.}\label{2circles-t}
\end{figure}

Let $X$ (see Fig.~\ref{2circles-d}) be a
point on a circle of $\P$ shown in solid line,
inside a disk $\D$ bounded by the circle $C$ of $\P$ shown in dashed line.
For $\mu\ge 0$, let $T_\mu$ (see Fig.~\ref{2circles-t}) be a primitive triangle having a non-integer corner $a_0$
mapping to $N$, a non-integer corner $a_1$ mapped to $S$ (resp., to $N$) when $\mu$ is even (resp., odd),
and an integer corner $a_2$ of order $\mu+1$ mapped to $X$.
The small black dots in Fig.~\ref{2circles-t} indicate the preimages of the vertices of $\P$
which are not corners of $T_\mu$ (though they are vertices of its net).
The angle at the corner $a_1$ of $T_\mu$ is equal (resp., complementary) to the angle at its corner $a_0$
when $\mu$ is even (resp., odd).
Then $T_\mu$ has the empty $(3,1)$-type
when $\mu=0$, the maximal $(3,1)$-type $\{(2,1)\}$ when $\mu>0$, and the multiplicity $\mu_{2,1}=\mu$.

\begin{cor}\label{ex:triangle}
Every primitive triangle with two non-integer corners over the partition $\P$ is combinatorially
equivalent to one of the triangles $T_\mu$.
A triangle $\bar T_\mu$ with the $(3,2)$-type can be obtained from the triangle $T_\mu$ by
reflection symmetry, relabeling the corners $a_1$ and $a_2$.
\end{cor}

\begin{figure}
\centering
\includegraphics[width=3.0in]{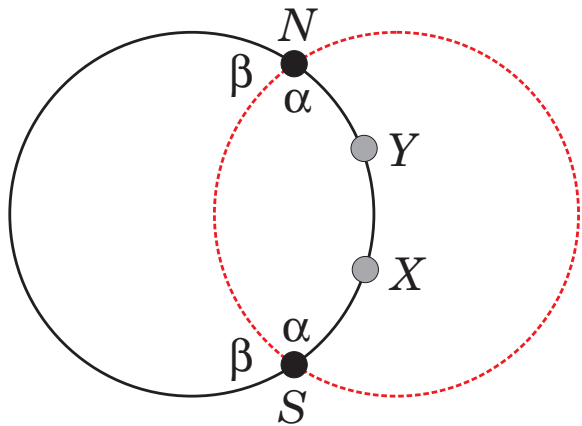}
\caption{Location of the points $X$ and $Y$ for adjacent integer corners.}\label{2circles-da}
\end{figure}

\begin{figure}
\centering
\includegraphics[width=5.0in]{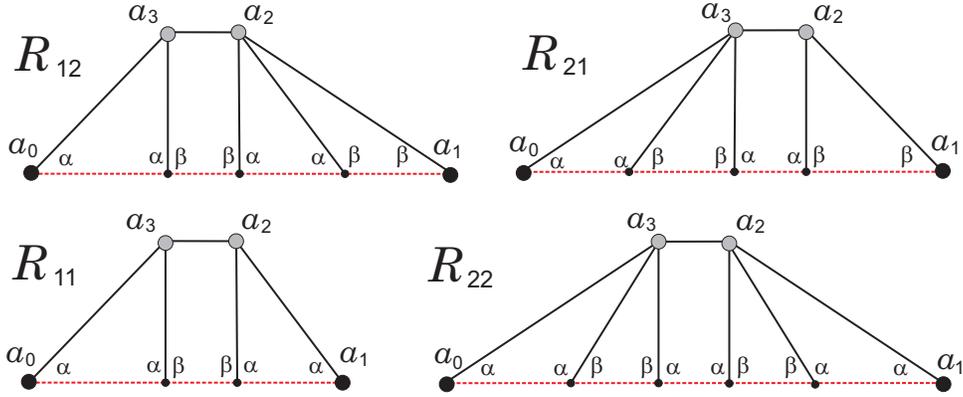}
\caption{Primitive quadrilaterals $R_{\mu\nu}$.}\label{2circles-r}
\end{figure}

Let $X$ and $Y$ (see Fig.~\ref{2circles-da}) be two points on the same arc of a circle of $\P$ shown in solid line,
inside a disk $\D$ bounded by the circle $C$ of $\P$ shown by the dashed line.
For $\mu,\nu\ge 0$, let $R_{\mu\nu}$ (see Fig.~\ref{2circles-r}) be a primitive quadrilateral having
a non-integer corner $a_0$ mapping to $N$, a non-integer corner $a_1$ mapped to $S$ (resp., to $N$)
when $\mu+\nu$ is even (resp., odd),
and integer corners $a_2$ of order $\nu+1$ and $a_3$ of order $\mu+1$ mapped to $X$ and $Y$
(resp., to $Y$ and $X$) when $\mu$ is even (resp., odd).
The side $L_1$ of $R_{\mu\nu}$ is mapped to the circle $C$ traversing it counterclockwise.
The angle at the corner $a_1$ of $R_{\mu\nu}$ is equal (resp., complementary) to the angle at its corner $a_0$
when $\mu+\nu$ is even (resp., odd).
Then $R_{\mu\nu}$ has the empty $(4,1)$-type when $\mu=\nu=0$, the $(4,1)$-type $\{(3,1)\}$ with the multiplicity
$\mu_{3,1}=\mu$ when $\mu>0,\,\nu=0$, the $(4,1)$-type $\{(2,1)\}$ with the multiplicity
$\mu_{2,1}=\nu$ when $\mu=0,\,\nu>0$, and the (unique) maximal $(4,1)$-type $\{(3,1),(2,1)\}$
with the multiplicities $\mu_{3,1}=\mu$ and $\mu_{2,1}=\nu$ when $\mu,\nu>0$.
Note that when $\mu=0$ (resp., $\nu=0$) the corner $a_3$ (resp., $a_2$) of $R_{\mu\nu}$
is removable, thus $R_{\mu\nu}$ is isometric to a triangle (to a digon when $\mu=\nu=0$).
We'll need such quadrilaterals later as building blocks for constructing non-primitive
quadrilaterals.

\begin{cor}\label{ex:quad-abcd}
Every primitive quadrilateral over $\P$ with two adjacent non-integer corners
is combinatorially equivalent to one of the quadrilaterals $R_{\mu\nu}$.
A quadrilateral $\bar R_{\mu\nu}$ with the $(4,3)$-type can be obtained from the quadrilateral $R_{\mu\nu}$
by reflection symmetry, relabeling the corners $a_1,a_2,a_3$.
\end{cor}

\begin{figure}
\centering
\includegraphics[width=4.0in]{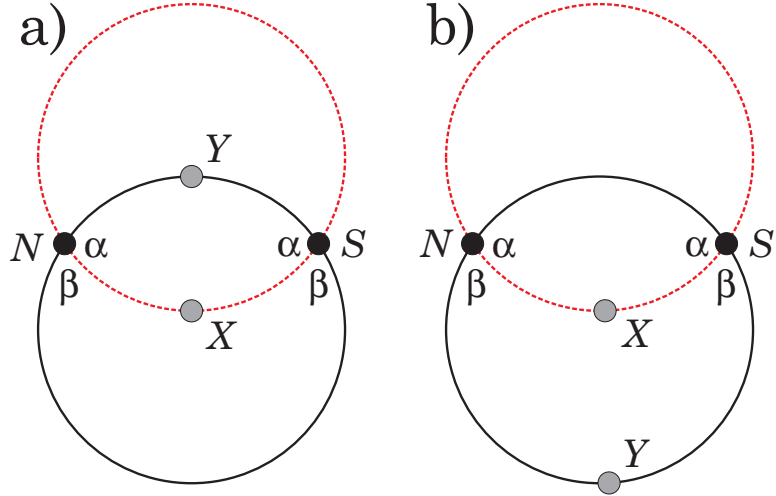}
\caption{Location of the points $X$ and $Y$ for two opposite integer corners.}\label{2circles2a}
\end{figure}

Let $X$ and $Y$ (see Fig.~\ref{2circles2a}) be two points on distinct circles of $\P$.

\begin{figure}
\centering
\includegraphics[width=5.0in]{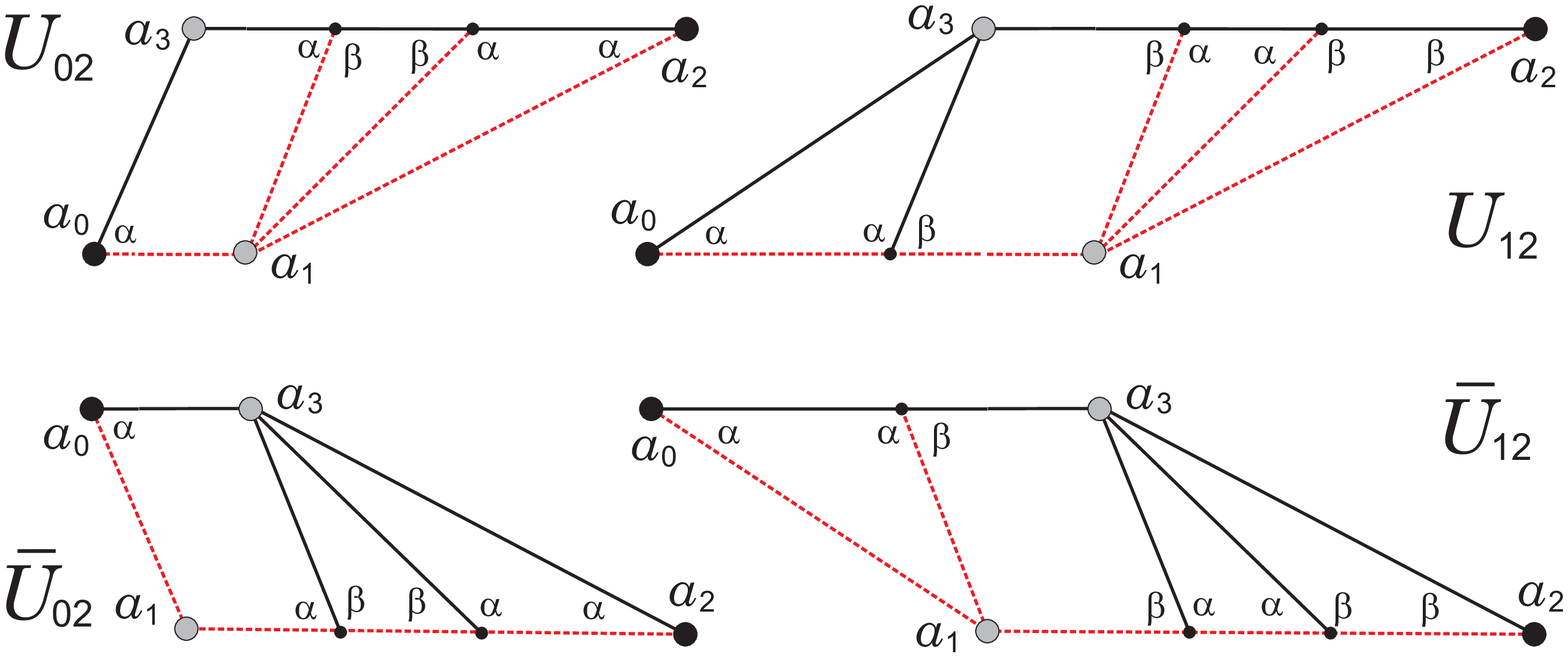}
\caption{Primitive quadrilaterals $U_{\mu\nu}$ and $\bar U_{\mu\nu}$.}\label{2circles-u1}
\end{figure}

For $\mu,\nu\ge 0$, let $U_{\mu\nu}$ and $\bar U_{\mu\nu}$ (see Fig.~\ref{2circles-u1})
be primitive quadrilaterals having a non-integer corner $a_0$
mapping to $N$, a non-integer corner $a_2$ mapped to $S$ (resp., to $N$) when $\mu+\nu$ is even (resp., odd),
and integer corners $a_1$ of order $\mu+1$ and $a_3$ of order $\nu+1$ mapped to $X$ and $Y$, respectively,
as shown in Fig.~\ref{2circles2a}a if $\mu$ is even and Fig.~\ref{2circles2a}b if $\mu$ is odd.
The angle at the corner $a_2$ of $U_{\mu\nu}$ and $\bar U_{\mu\nu}$ is equal (resp., complementary)
to the angle at its corner $a_0$ when $\mu+\nu$ is even (resp., odd).
Then $U_{\mu\nu}$ (resp., $\bar U_{\mu\nu}$) has the empty $(4,2)$-type when $\mu=\nu=0$,
the $(4,2)$-type $\{(3,1)\}$ (resp., $\{(1,4)\}$) with the multiplicity
$\mu_{3,1}=\mu$ (resp., $\mu_{1,4}=\mu$) when $\mu>0,\,\nu=0$, the $(4,2)$-type $\{(1,3)\}$
(resp., $\{(3,2)\}$) with the multiplicity
$\mu_{1,3}=\nu$ (resp., $\mu_{3,2}=\nu$) when $\mu=0,\,\nu>0$, and the maximal $(4,2)$-type $\{(3,1),(1,3)\}$
(resp., $\{(1,4),(3,2)\}$) with the multiplicities $\mu_{3,1}=\mu$ and $\mu_{1,3}=\nu$
(resp., $\mu_{1,4}=\mu$ and $\mu_{3,2}=\nu$) when $\mu,\nu>0$.
Note that when either $\mu=0$ or $\nu=0$, the quadrilateral $U_{\mu\nu}$ and $\bar U_{\mu\nu}$
has a removable integer corner (both integer corners when $\mu=\nu=0$)
and is isometric to a triangle (a digon when $\mu=\nu=0$).
We'll need such quadrilaterals later as building blocks for constructing non-primitive
quadrilaterals.

\begin{figure}
\centering
\includegraphics[width=5.0in]{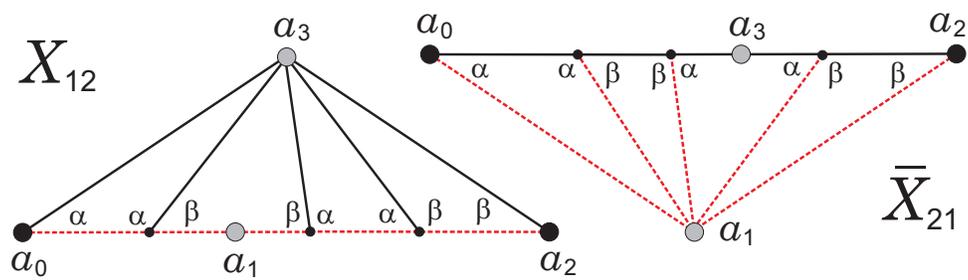}
\caption{Primitive quadrilaterals $X_{\mu\nu}$ and $\bar X_{\mu\nu}$.}\label{2circles-x}
\end{figure}

For $\mu,\nu\ge 0$, let $X_{\mu\nu}$ and $\bar X_{\mu\nu}$ (see Fig.~\ref{2circles-x})
be primitive quadrilaterals having a non-integer corner $a_0$
mapping to $N$, a non-integer corner $a_2$ mapped to $S$ (resp., to $N$) when $\mu+\nu$ is even (resp., odd),
and integer corners $a_1$ of order $1$ and $a_3$ of order $\mu+\nu+1$ mapped to $X$ and $Y$, respectively,
as shown in Fig.~\ref{2circles2a}a if $\mu$ is even and Fig.~\ref{2circles2a}b if $\mu$ is odd.
The angle at the corner $a_2$ of $X_{\mu\nu}$ and $\bar X_{\mu\nu}$ is equal (resp., complementary)
to the angle at its corner $a_0$ when $\mu+\nu$ is even (resp., odd).
Then $X_{\mu\nu}$ (resp., $\bar X_{\mu\nu}$) has the empty $(4,2)$-type when $\mu=\nu=0$,
the $(4,2)$-type $\{(3,1)\}$ (resp., $\{(1,4)\}$) with the multiplicity
$\mu_{3,1}=\mu$ (resp., $\mu_{1,4}=\mu$) when $\mu>0,\,\nu=0$, the $(4,2)$-type $\{(3,2)\}$
(resp., $\{(1,3)\}$) with the multiplicity
$\mu_{3,2}=\nu$ (resp., $\mu_{1,3}=\nu$) when $\mu=0,\,\nu>0$, and the maximal $(4,2)$-type $\{(3,1),(3,2)\}$
(resp., $\{(1,4),(1,3)\}$) with the multiplicities $\mu_{3,1}=\mu$ and $\mu_{3,2}=\nu$
(resp., $\mu_{1,4}=\mu$ and $\mu_{1,3}=\nu$) when $\mu,\nu>0$.
Note that the corner $a_1$ of $X_{\mu\nu}$ and the corner $a_3$ of $\bar X_{\mu\nu}$
are removable, thus these quadrilaterals are isometric to triangles $T_{\mu+\nu}$.
We'll need such quadrilaterals later as building blocks for constructing non-primitive quadrilaterals.

\begin{cor}\label{ex:quad-acbd}
Every primitive quadrilateral over $\P$ with two opposite non-integer corners
is combinatorially equivalent to one of the quadrilaterals $U_{\mu\nu}$,
$\bar U_{\mu\nu}$, $X_{\mu\nu}$, $\bar X_{\mu\nu}$.
\end{cor}

Note that $U_{\mu0}$ and $X_{\mu0}$ are combinatorially equivalent for each $\mu$,
$U_{0\nu}$ and $\bar X_{0\nu}$ are combinatorially equivalent for each $\nu$,
$\bar U_{\mu0}$ and $\bar X_{\mu0}$ are combinatorially equivalent for all $\mu$,
$\bar U_{0\nu}$ and $X_{0,\nu}$ are combinatorially equivalent for all $\nu$.

\begin{prop}\label{primitive-count}
For given $n\ge 2$ and $k,\;0<k<n$, the number $M(n,k)$ of distinct maximal
$(n,k)$-types satisfies the following recurrence:
\begin{equation}\label{max-recurrence}
M(n,k)=\sum_{m=1}^k M(n-m,k-m+1)+\sum_{m=1}^{n-k} M(n-m,n-k-m+1).
\end{equation}
Since $M(2,1)=1$, this implies that
\begin{equation}\label{max-generating}
\sum_{k,l=1}^\infty M(k+l,k)x^k y^l=\frac{(1-x)(1-y)}{(1-x)(1-y)-x(1-y)-y(1-x)}.
\end{equation}
\end{prop}

\medskip
\noindent{\bf Proof.} This recurrence follows from construction in the
proof of Lemma \ref{maxtype}.

\section{Irreducible spherical polygons with two\newline non-integer corners}\label{irreducible-2}
In this section, $Q$ is an irreducible marked $n$-gon with two non-integer corners.
Its corners and sides are labeled as in the previous section.
In particular, the non-integer corners of $Q$ are $a_0$ and $a_k$.
Let $\Q$ be the net of $Q$, and let $\Q_C$ and $\Q_{C'}$ be preimages
of the circles $C$ and $C'$ of $\P$, respectively.
Since $Q$ is irreducible, its net $\Q$ does not contain an
interior arc with the ends at two distinct corners of $Q$.
However, $\Q$ may have loops.

\begin{lemma}\label{corners}
Integer corners of $Q$ do not belong to preimages of the vertices of $\P$.
\end{lemma}

\medskip{\bf Proof.}
We proceed by induction on the number $m$ of faces of $\Q$, the case $m=1$ being trivial.
Let $p$ be an integer corner of $Q$ that belongs to the preimage of a vertex of $\P$.
Then there is a maximal interior arc $\gamma$ of $\Q$ with an end at $p$,
such that $\gamma$ partitions the angle at $p$ into two non-integer angles.
The other end $q$ of $\gamma$ cannot be on a side of $Q$.
Otherwise, $\gamma$ would partition $Q$ into two irreducible
polygons, each having non-integer angles at both $p$ and $q$,
which contradicts Theorem \ref{non-integer}.
Thus $q=p$ and $\gamma$ is a loop.
Let $\gamma$ be an arc of $\Q_C$.
Then there is a maximal interior arc $\gamma'$ of $\Q_{C'}$
having one end at $p$ and intersecting $\gamma$ in an interior vertex $q$ of $\Q$.
We saw in the proof of Theorem \ref{non-integer} that $\gamma'$ cannot be a loop.
Thus $\gamma'$ partitions $Q$ into two irreducible polygons, each of them
having less than $m$ faces in its net,
with $p$ being an integer corner of each of them.
This contradicts the induction hypothesis.

\begin{lemma}\label{loop}
Let $\gamma$ be a loop of $\Q$.
Then $\gamma$ contains an integer corner of $Q$.
\end{lemma}

\medskip
\noindent{\bf Proof.}
Let $\gamma$ be a loop in $\Q_C$ that
does not contain a corner of $Q$.
Then $\gamma$ intersects $\Q_{C'}$ at two points
$q$ and $q'$. Let $\gamma'$ be the maximal loop in $\Q_{C'}$
passing through $q$ and $q'$.
The same arguments as in the proof of Theorem \ref{non-integer}
show that $\gamma'$ cannot be a loop.
Since $Q$ is irreducible, $\gamma'$ cannot have both ends
at the corners of $Q$. Thus one of its ends is at the side of $Q$.
But this implies that a face of $\Q$ adjacent to $\gamma$
outside the disk bounded by $\gamma$ has a disconnected
intersection with $\Q_C$, a contradiction.

\begin{lemma}\label{twoloops}
Let $\gamma$ be a loop of $\Q_C$ with both ends at an integer corner $p$ of $Q$.
Then $\gamma$ there is a loop $\gamma'$
of $\Q_{C'}$ intersecting $\gamma$ at two points and having both ends at an
integer corner $p'$ of $Q$.
The union of $\gamma$ and $\gamma'$ is a
pseudo-diagonal of $Q$ shown in Fig.~\ref{2circles-i}.
\end{lemma}

\begin{df}\label{pseudo}
The number $\nu_{ij}$ of the pseudo-diagonals of $Q$
(see Fig.~\ref{2circles-i}) connecting integer corners $a_i$ and $a_j$ of $Q$
is called the {\em multiplicity} of a pseudo-diagonal
connecting $a_i$ and $a_j$.
\end{df}

\begin{thm}\label{irreducible}
Each irreducible spherical polygon $Q$ with two non-integer corners
can be obtained from a primitive polygon by adding
(multiple) pseudo-diagonals connecting some of its integer corners.
These pseudo-diagonals do not intersect inside $Q$.
For a primitive polygon with a maximal type,
the irreducible polygons that can be obtained from it
are uniquely determined by the multiplicities of
the pseudo-diagonals.
\end{thm}

\medskip
\noindent{\bf Proof.}
The first part of the statement is obvious, since removing
all loops from the net of $Q$ we obtain a primitive polygon $Q'$.
If $Q'$ has maximal type then each face of its net $\Q'$
has at most two integer corners of $Q'$ in its boundary.
The multiplicities of the pseudo-diagonals connecting
these pairs of corners completely determine the irreducible
polygon $Q$ from which $Q'$ was obtained.

\medskip
Each of the primitive quadrilaterals $U_{\mu\nu}$, $\bar U_{\mu\nu}$,
$X_{\mu\nu}$, $\bar X_{\mu\nu}$ (see Figs.~\ref{2circles-u1},\ref{2circles-x})
contains a single face $F$ of its net with both integer corners $a_1$ and $a_3$
in its boundary. Adding $\kappa$ pseudo-diagonals connecting the two
integer corners inside $F$, we obtain irreducible quadrilaterals
$U_{\mu\nu}^\kappa$, $\bar U_{\mu\nu}^\kappa$, $X_{\mu\nu}^\kappa$, $\bar X_{\mu\nu}^\kappa$.
For $\kappa>0$, these quadrilaterals are not primitive.
We identify $U_{\mu\nu}^0$, $\bar U_{\mu\nu}^0$, $X_{\mu\nu}^0$, $\bar X_{\mu\nu}^0$ with
$U_{\mu\nu}$, $\bar U_{\mu\nu}$, $X_{\mu\nu}$, $\bar X_{\mu\nu}$, respectively.

\begin{cor}\label{ex:irreducible}
Every marked irreducible quadrilateral over $\P$ with two opposite non-integer
corners is combinatorially equivalent to one of the quadrilaterals
$U_{\mu\nu}^\kappa$, $\bar U_{\mu\nu}^\kappa$, $X_{\mu\nu}^\kappa$, $\bar X_{\mu\nu}^\kappa$
with $\mu,\nu,\kappa\ge 0$.
\end{cor}

Note that $U_{\mu0}^\kappa$ and $X_{\mu0}^\kappa$ are combinatorially equivalent for each $\mu,\kappa$,
$U_{0\nu}^\kappa$ and $\bar X_{0\nu}^\kappa$ are combinatorially equivalent for each $\nu,\kappa$,
$\bar U_{\mu0}^\kappa$ and $\bar X_{\mu0}^\kappa$ are combinatorially equivalent for all $\mu,\kappa$,
$\bar U_{0\nu}^\kappa$ and $X_{0,\nu}^\kappa$ are combinatorially equivalent for all $\nu,\kappa$.

\begin{df} The $(n,k)$-{\em type} of an irreducible polygon $Q$
(or simply the type of $Q$ if $n$ and $k$ are not specified)
is the type of the primitive polygon $Q'$
obtained from $Q$ by removing all loops,
with the unordered pairs $(i,j)$ added for the
integer corners $a_i$ and $a_j$ with positive
multiplicities $\nu_{ij}$ of the pseudo-diagonals
connecting these corners.
The type of $Q$ is {\em maximal} if the type of $Q'$
is maximal, and any two integer corners of $Q'$
belonging to the boundary of the same face of its net
are connected in $Q$ by at least one pseudo-diagonal.
Every $(n,k)$-type of an irreducible spherical polygon can be obtained
from a (non-unique) maximal $(n,k)$-type if some of
the multiplicities are allowed to be zero.
\end{df}

\begin{prop}
The number of distinct maximal $(n,k)$ types of irreducible spherical
polygons equals the number of distinct maximal $(n,k)$ types of primitive
spherical polygons.
\end{prop}

\medskip
\noindent{\bf Proof.} This follows from Theorem \ref{irreducible}.

\section{Classification of spherical digons and\newline triangles}\label{sec:classification}

For $m\ge 1$, there is a unique rational spherical digon $D_m$, each of its corners equal $m\pi$.

There is a unique, up to combinatorial equivalence, irreducible spherical
digon with two equal non-integer corners, isometric to a face of $\P$.
Any irrational spherical digon is obtained from it by attaching a digon $D_m$.
Here $m\ge 0$, with $m=0$ meaning that nothing is attached.

Theorem \ref{irreducible} implies that any irreducible spherical triangle with
two non-integer corners is primitive and
combinatorially equivalent to one of the triangles $T_\mu$ (see Corollary \ref{ex:triangle}).

Any spherical triangle $T$ over $\P$ is either rational or has two non-integer corners.
If $T$ has two non-integer corners $a_0$ and $a_1$, it is combinatorially equivalent to a
triangle $T_\mu$ with digons $D_i$, $D_j$ and $D_l$ attached to its sides
$L_1$, $L_2$ and $L_3$, respectively, where $i,j,l\ge 0$, the value $0$ meaning that
no digon is attached, and $i>0$ only when $\mu=0$.

\begin{figure}
\centering
\includegraphics[width=2.5in]{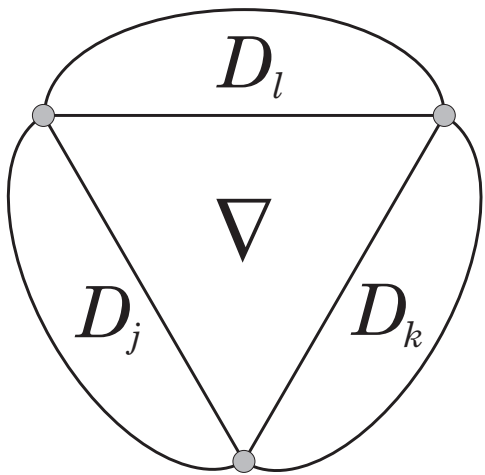}
\caption{Rational triangle $\nabla_{jkl}$.}\label{nabla}
\end{figure}

If $T$ is rational, it is combinatorially equivalent
to a rational triangle $\nabla$ with three removable corners, with digons $D_j$, $D_k$ and $D_l$
attached to its three sides, where $j,k,l\ge 0$ are determined by $T$ uniquely up to cyclic permutation.
We use notation $\nabla_{jkl}$ for such a triangle (see Fig.~\ref{nabla}).

\section{Classification of spherical quadrilaterals\newline
with two adjacent non-integer corners}\label{classification-adjacent}

We consider marked quadrilaterals over $\P$ with two adjacent non-integer corners $a_0$ and $a_1$.

\begin{lemma}
A quadrilateral $Q$ over $\P$ cannot have more than two non-integer
corners.
\end{lemma}

\medskip
\noindent{\bf Proof.}
If $Q$ is a union of irreducible polygons, then one of them, say $\tilde Q$, should be either
a quadrilateral or a triangle.
If $\tilde Q$ is a quadrilateral,
then, by Lemma \ref{corners}, $\tilde Q$ must have two corners that are not mapped to
vertices of $\P$. These corners of $\tilde Q$, and the corresponding corners of $Q$, must be integer.
If $\tilde Q$ is a triangle, at least one of its corners is not mapped to
a vertex of $\P$, thus $Q$ has at least one integer corner.
Since the number of non-integer corners of any polygon over $\P$ is even,
at least two corners of $Q$ are integer.

\medskip
A quadrilateral $Q$ can be partitioned into irreducible polygons.
We'll see later that for a quadrilateral with two adjacent non-integer corners
this partition is unique.

\begin{figure}
\centering
\includegraphics[width=5.0in]{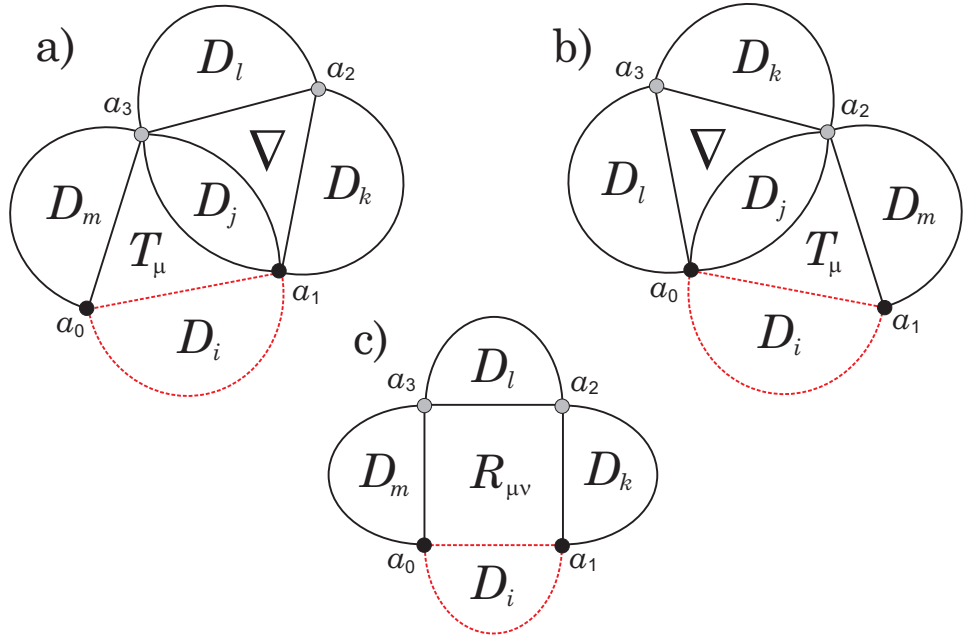}
\caption{Quadrilaterals with adjacent non-integer corners.}\label{2circles-abcd}
\end{figure}

If one of these polygons is a triangle $T_\mu$ with its corners
at $a_0$, $a_1$ and either $a_2$ or $a_3$,
then $Q$ is a union of $T_\mu$, a triangle $\nabla_{jkl}$ attached to its side other than $L_1$
so that $D_j$ is adjacent to $T_\mu$,
and digons $D_i$ and $D_m$ attached to the other two sides of $T_\mu$
so that  $D_i$ is adjacent to its side $L_1$ (see Fig.~\ref{2circles-abcd}ab).
Here $i,j,k,l,m\ge 0$, with $i>0$ only if $\mu=0$.

Otherwise, $Q$ contains a quadrilateral
$R_{\mu\nu}$ having the same corners as $Q$.
Then $Q$ is the union of $R_{\mu\nu}$ and digons $D_i$, $D_k$, $D_l$, $D_m$
attached to its sides (see Fig.~\ref{2circles-abcd}c).
Here $i,j,k,l\ge 0$, with $l>0$ only if $\mu=\nu=0$.

\begin{figure}
\centering
\includegraphics[width=3.0in]{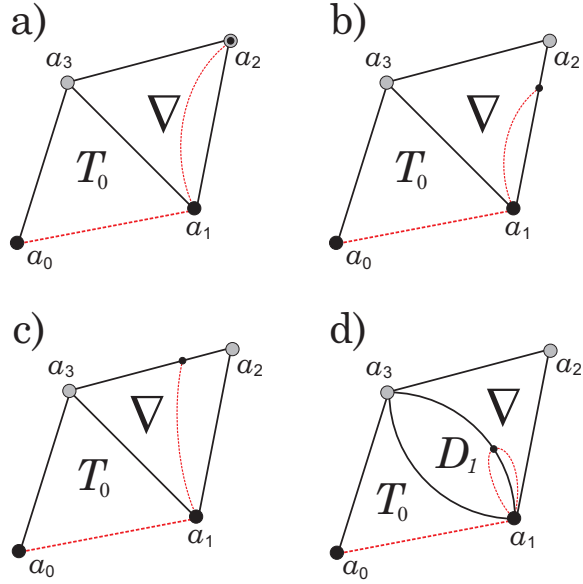}
\caption{Quadrilaterals in Fig.~\ref{2circles-abcd}a with $j=0$ (a-c) and $j=1$ (d).}\label{2circles-abcd-tnabla}
\end{figure}

Note that a vertex of $\nabla$ that is not a vertex of $T_{\mu}$
may be mapped to a vertex of $\P$ when $j$ is even, but not when $j$ is odd.
Fig.~\ref{2circles-abcd-tnabla}a-c shows the complete net for the quadrilateral
in Fig.~\ref{2circles-abcd}a with $j=0$ (and $\mu=i=k=l=m=0$).
All these quadrilaterals are combinatorially equivalent, although the corner
$a_2$ is mapped to a vertex of $\P$ in Fig.~\ref{2circles-abcd-tnabla}a
but not in Fig.~\ref{2circles-abcd-tnabla}bc.
Fig.~\ref{2circles-abcd-tnabla}d shows the complete net for
the quadrilateral in Fig.~\ref{2circles-abcd}a with $j=1$.
In this case, the corner $a_2$ cannot be mapped to a vertex of $\P$,
since there are two points mapped to vertices of $\P$ on
the side of $\nabla$ connecting $a_1$ and $a_3$.

\medskip{\bf Counting quadrilaterals with given angles.}
We want to classify marked spherical quadrilaterals $Q$
with two adjacent non-integer corners $a_0$ and $a_1$,
and with given orders $A_0,\dots,A_3$ of its corners,
up to combinatorial equivalence.
Here $A_0,A_1\ge 0$ and $A_2,A_3>0$ are the orders of non-integer
and integer corners of $Q$, respectively.

We define the following numbers: $\delta=\frac12(A_1+A_3-A_0-A_2)$,
$\sigma=\frac12(A_2+A_3-A_0-A_1)$, $\rho=\frac12(A_2+A_3-|A_1-A_0|)$.
These numbers are integer if and only if the
corners $a_0$ and $a_1$ of $Q$ are mapped to distinct vertices of $\P$.

The relations of these parameters with those in Theorem \ref{theorem3} is
the following
$$A_0=[\alpha_3],\quad A_1=[\alpha_0],\quad A_2=\alpha_1,\quad A_3=\alpha_2,$$
In Case a), $\rho=\kappa+1$, and in Case b) $\sigma=\kappa+3/2.$

\begin{lemma}\label{count-lemma}
For given positive integers $p$, $q$, $r$, $s$ satisfying $p+r=q+s$,
the system of equations $x+y=p+1,\;y+z=q+1,\;z+t=r+1,\;t+x=s+1$ has
$\min(p,q,r,s)$ solutions $(x,y,z,t)$ in positive integers.
\end{lemma}

The proof is left as an easy exercise for the reader.

\begin{prop}\label{count-abcd}
A marked quadrilateral $Q$ with non-integer corners $a_0$ and $a_1$
mapped to distinct vertices of $\P$ exists if and only if $\rho$ is a positive integer.
In this case, there are $\min(A_2,A_3,\rho)$ combinatorially distinct quadrilaterals with given
orders $A_0,\dots,A_3$ of their angles.

A quadrilateral $Q$ with the corners $a_0$ and $a_1$ mapped to the same vertex of $\P$
exists if and only if $\sigma-1$ is a positive non-integer.
In this case, there are $\min(A_2,A_3,[\sigma])$ combinatorially distinct quadrilaterals with given angles.
\end{prop}

\medskip
\noindent{\bf Proof.}
For a quadrilateral $Q$ shown in Fig.~\ref{2circles-abcd}a,
its corners $a_0$ and $a_1$ are mapped to distinct vertices of $\P$ if $\mu$ is even,
and to the same vertex of $\P$ if $\mu$ is odd. Orders of the corners of $Q$
are $A_0=i+m$, $A_1=i+j+k+1$, $A_2=k+l+1$, $A_3=j+k+l+m+1$.
We have $\delta=j+1+\mu/2$, $\sigma=l-i+1+\mu/2$.
If $A_1>A_0$, that is, $j+k\ge m$, then $\rho=l+m+1+\mu/2$,
otherwise, $\rho=j+k+l+2+\mu/2>A_2$.
In particular, $\delta$ is positive.
Similarly, $\delta$ is negative for the quadrilateral $\Q$ shown in Fig.~\ref{2circles-abcd}b.

For the quadrilateral $Q$ in Fig.~\ref{2circles-abcd}c, its corners
$A$ and $B$ are mapped to distinct vertices of $\P$ if $\mu+\nu$ is even,
and to the same vertex of $\P$ if $\mu+\nu$ is odd.
Orders of the corners of $Q$ are $A_0=i+m$, $A_1=i+k$, $A_2=k+l+1+\nu$, $A_3=l+m+1+\mu$.
We have $\delta=(\mu-\nu)/2$, $\sigma=l-i+1+(\mu+\nu)/2$.
If $A_1\ge A_0$, that is, $k\ge m$, then $\rho=l+m+1+(\mu+\nu)/2$,
otherwise $\rho=k+l+1+(\mu+\nu)/2$.

Note that in all cases $\delta$, $\sigma$ and $\rho$ are integer if and only if
the corners $a_0$ and $a_1$ are mapped to distinct vertices of $\P$.
If $i>0$ then $\mu=0$ in Fig.~\ref{2circles-abcd}ab and $\mu=\nu=0$ in Fig.~\ref{2circles-abcd}c,
thus the fractional parts of the angles at $a_0$ and $a_1$ are equal.
If $i=0$ then $\sigma\ge 1$.

We start with the quadrilaterals with equal fractional parts of the angles at $a_0$ and $a_1$.

Consider first the case $\delta=0$, that is, $A_3-A_0=A_2-A_1$.
This is only possible for a quadrilateral $\Q$ in Fig.~\ref{2circles-abcd}c with $\mu=\nu$.
We may assume $A_0-A_1=A_3-A_2\ge 0$ (the other case follows by symmetry).
Then $A_3\ge A_2$ and $\rho=(A_3+A_2-A_0+A_1)/2=A_2$, so we have to prove that the number of
combinatorially distinct quadrilaterals is $A_2$.

Lemma \ref{count-lemma} applied to $x=i$, $y=m+1$,
$z=l+1$, $t=k+1$ implies that the number of quadrilaterals with $i>0$
is $\min(A_0,A_1,A_2,A_3)=\min(A_1,A_2)$.
If $A_1\ge A_2$, the number of quadrilaterals is $A_2$.
Otherwise, there are $A_2-A_1=A_3-A_0$ quadrilaterals with $i=0$, $m=A_0$, $k=A_1$,
$0\le \mu=\nu<A_2-A_1$, $l=A_2-A_1-1-\nu=A_3-A_0-1-\mu$.
Thus the total number of quadrilaterals is again $C$.

Consider now the case $\delta>0$ (the case $\delta<0$ follows by symmetry).
Then $A_3>A_2+A_0-A_1$, so $A_3>A_2$ and $\rho=(A_2+A_3+A_1-A_0)/2=A_2+\delta>A_2$ if $A_0\ge A_1$.
If $A_0<A_1$ then $\rho=(A_2+A_3+A_0-A_1)/2=A_3-\delta<A_3$.
Thus we have to prove that the number of quadrilaterals is $\min(A_2,A_3-\delta)$.

The following three subcases are possible:\newline
(i) a quadrilateral in Fig.~\ref{2circles-abcd}c with $i=0$ and $\mu>\nu$,\newline
(ii) a quadrilateral in Fig.~\ref{2circles-abcd}a with $i=0$,\newline
(iii) a quadrilateral in Fig.~\ref{2circles-abcd}a with $i>0$ and $\nu=0$.

Subcase (i). For a quadrilateral $Q$ in Fig.~\ref{2circles-abcd}c with $i=0$ and $\mu>\nu$,
we have $m=A_0$, $k=A_1$, $A_2>A_1$.
For each $\nu$ such that $0\le\nu<A_2-A_1$, there is a quadrilateral with $l=A_2-A_1-1-\nu$ and $\mu=\nu+2\delta$.
Thus there are $A_2-A_1$ quadrilaterals in this case.

Subcase (ii). For a quadrilateral $Q$ in Fig.~\ref{2circles-abcd}a with $i=0$, we have
$m=A_0$, $j+k+1=A_1$, $k+l+1=A_2$, $\mu+j+l+2=A_3-A_0$, $\delta=j+1+\mu/2$, $\sigma=l+1+\mu/2>0$.
For a fixed even $\mu\ge 0$, we have
$j=\delta-1-\mu/2$, $l=\sigma-1-\mu/2$, $k=A_1-\delta+\mu/2$.
Thus $\delta-A_1\le\mu/2\le\min(\delta,\sigma)$.
Since $\mu$ is even, the number of quadrilaterals in this case is
$\min(A_1,A_2)$ if $A_1\le\delta$ and $\min(\delta,\sigma)$ if $A_1>\delta$.
Note that $\sigma-\delta=A_2-A_1$.

Subcase (iii). For a quadrilateral $Q$ in Fig.~\ref{2circles-abcd}a with $i>0$ and $\mu=0$, we have
$\delta=j+1$, $i+m=A_0$, $i+k=A_1-\delta$, $k+l+1=A_2$, $l+m+1=A_3-\delta$.
Lemma \ref{count-lemma} applied to $x=i$, $y=m+1$, $z=l+1$, $t=k+1$ implies that the number
of quadrilaterals is $\min(A_0,A_1-\delta,A_2,A_3-\delta)$.
If $A_2>A_1$, that is $l>i+j$, then $A_3-\delta>A_0$ and the number of quadrilaterals is $\min(A_0,A_1-\delta)$.

Combining subcases (i)-(iii), we see first that, if $A_2>A_1$ and $A_1\le\delta$,
then the number of quadrilaterals is $A_2$ (there are no quadrilaterals in the subcase (iii) in this case).
Note that $A_3-\delta=(A_3+A_2+A_0-A_1)/2=A_2+\delta-A_1\ge A_2$ in this case.

Next, if $A_2>A_1$ and $A_1>\delta$, the total number of quadrilaterals is
$\min(A_2-A_1+\delta+A_0,A_2)=\min(A_3-\delta,A_2)$.

Finally, if $A_2\le A_1$ then there are no quadrilaterals in the subcase (i).
If $A_1\le\delta$, there are no quadrilaterals in the subcase (iii),
and the number of quadrilaterals in the subcase (ii) is $A_2$.
If $A_1>\delta$ and $\sigma>0$ then the number of quadrilaterals in the subcase (ii) is $\sigma$.
Since $A_0+\sigma=A_3-\delta$ and $A_1-\delta+\sigma=A_2$
the total number of quadrilaterals is $\min(A_2,A_3-\delta)$.
If $\sigma\le 0$, there are no quadrilaterals in subcases (i) and (ii),
$A_0\ge A_2+A_3-A_1\ge A_3$ and $A_1-\delta\ge A_2$.
Thus the number of quadrilaterals in (iii) is $\min(A_2,A_3-\delta)$.

If the fractional parts of the angles at $a_0$ and $a_1$ are complementary then $i=0$.
Since $\delta\ne 0$ in this case, we may assume $\delta>0$ (the other case follows by symmetry).
Since $i=0$, only subcases (i) and (ii) above are possible.

Repeating the above arguments, we see that $A_2>A_1$ in the subcase (i), and the number of quadrilaterals is $A_2-A_1$.
In the subcase (ii), since $\mu\ge 1$ is odd, the number of quadrilaterals is
$\min(A_1,A_2)$ when $A_1<\delta$ and $\min([\delta],[\sigma])$ when $A_1>\delta$.

Thus the total number of quadrilaterals is $A_2$ if $A_1<\delta$.
Since $A_2-A_1=\sigma-\delta$, we have $\sigma>A_2$ in this case.

When $A_1>\delta$, either $A_2\le A_1$ and the total number of quadrilaterals is $[\sigma]$
(since there are no quadrilaterals in the subcase (i))
or $A_2>A_1$ and the total number of quadrilaterals is again $[\sigma]=A_2-A_1+\delta$.
Note that $\sigma<A_2$ when $A_1>\delta$.

This completes the proof.

\section{Classification of spherical quadrilaterals\newline
with two opposite non-integer corners}\label{classification-opposite}

We consider now marked quadrilaterals $Q$ with two opposite
non-integer corners $a_0$ and $a_2$.
We assume that the corner $a_0$ is mapped to the vertex $N$ of $\P$.
The corner $a_2$ may be mapped either to the vertex $S$ or to the vertex $N$,
depending on the net of $Q$ (see Corollary \ref{ex:quad-acbd}).
For an irreducible quadrilateral $Q$, the integer corners $a_1$ and $a_3$
are mapped to the points $X$ and $Y$ on two distinct circles of $\P$
as shown on Fig.~\ref{2circles2a}.
However, for a reducible quadrilateral $Q$ one of these corners may be
mapped to a vertex of $\P$.
Note that a partition of such a quadrilateral $Q$ into
irreducible polygons may be non-unique.

\begin{figure}
\centering
\includegraphics[width=3.0in]{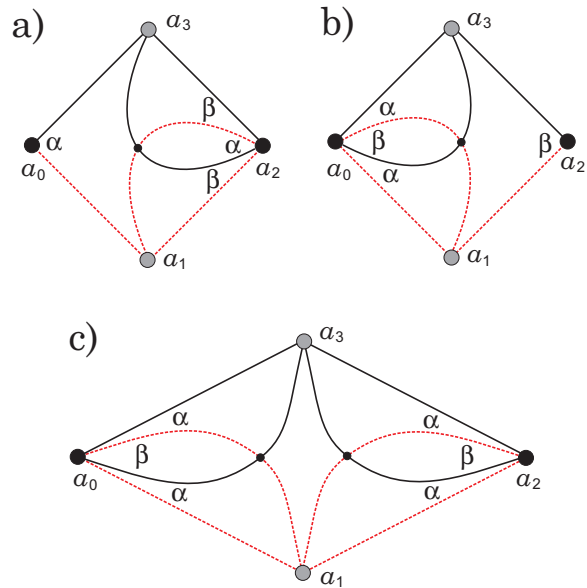}
\caption{Non-unique partition into irreducible polygons.}\label{fig:non-unique}
\end{figure}

\begin{example}\label{ex:non-unique}
{\rm Consider the quadrilaterals $I,J,K$ shown in Fig.~\ref{fig:non-unique}abc.
The quadrilateral $I$ can be represented in two different ways as the union of
a primitive quadrilateral (either $X_{01}$ or $\bar X_{01}$) and a
digon $D_1$ attached to its side of order 2.
The quadrilateral $J$ can be represented in two different ways as the union of
a primitive quadrilateral (either $X_{10}$ or $\bar X_{10}$) and a
digon $D_1$ attached to its side of order 2.
The quadrilateral $K$ can be represented in two different ways as the union of
a primitive quadrilateral (either $X_{22}$ or $\bar X_{22}$) and two digons $D_1$
attached to its adjacent sides of order 2.
Alternatively, $K$ can be represented in two different ways as the union
of a primitive quadrilateral (either $U_{22}$ or $\bar U_{22}$) and two digons $D_1$
attached to its opposite sides of order 2.}
\end{example}

\begin{figure}
\centering
\includegraphics[width=4.0in]{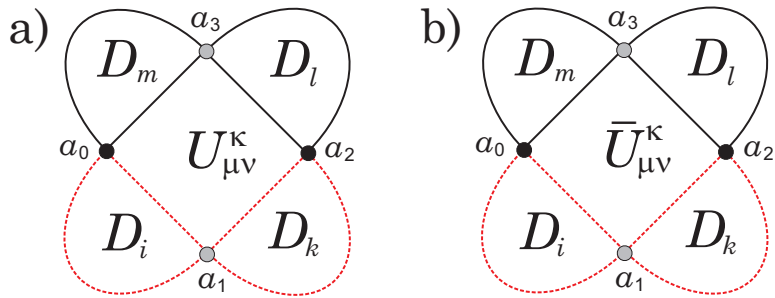}
\caption{Quadrilaterals of the types $U$ and $\bar U$.}\label{fig:acbd0}
\end{figure}

\begin{figure}
\centering
\includegraphics[width=4.0in]{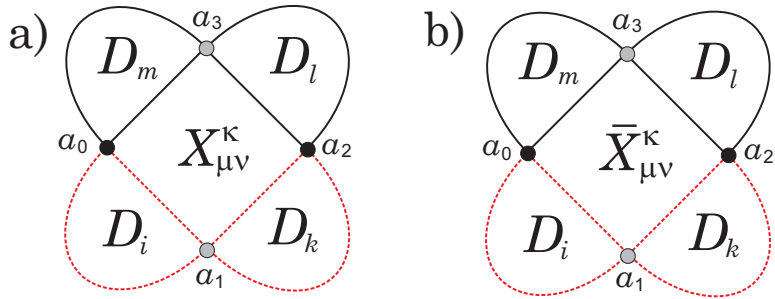}
\caption{Quadrilaterals of the types $X$ and $\bar X$.}\label{fig:acbd1}
\end{figure}

\begin{figure}
\centering
\includegraphics[width=4.0in]{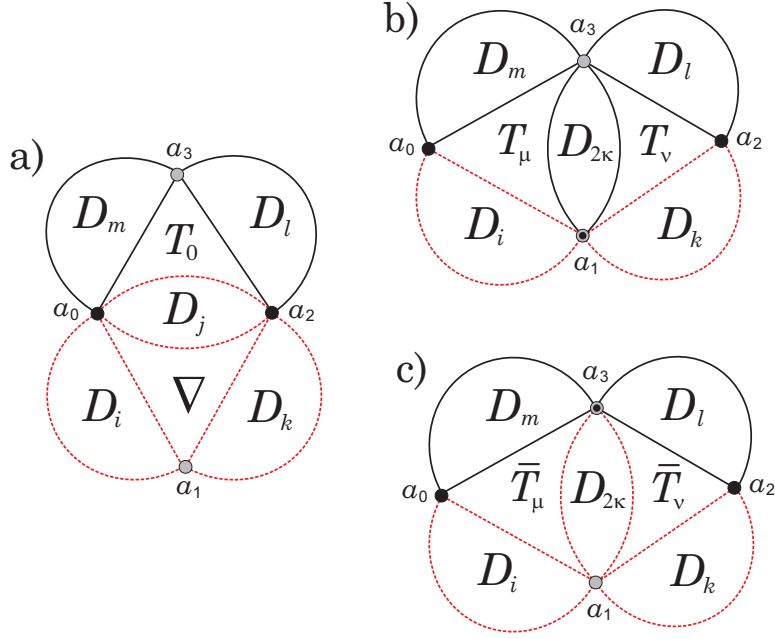}
\caption{Quadrilaterals of the types $T\nabla$, $TT$ and $\overline{TT}$.}\label{fig:acbd2}
\end{figure}
\newpage
\begin{prop}\label{quad-acbd}
Every marked spherical quadrilateral $Q$ over partition $\P$ with two opposite non-integer corners
is combinatorially equivalent to one of the following:

\noindent{\bf Type $U$.} An irreducible quadrilateral $U^\kappa_{\mu\nu}$ with digons $D_i$, $D_k$, $D_l$, $D_m$ attached to its sides
$L_1$, $L_2$, $L_3$, $L_4$ (see Fig.~\ref{fig:acbd0}a). Here $i,k,l,m,\mu,\nu,\kappa\ge 0$, with
$i>0$ only if $\mu\le 1$ and $l>0$ only if $\nu\le 1$.

\noindent{\bf Type $\bar U$.} An irreducible quadrilateral $\bar U^\kappa_{\mu\nu}$ with digons $D_i$, $D_k$, $D_l$, $D_m$ attached to its sides
$L_1$, $L_2$, $L_3$, $L_4$ (see Fig.~\ref{fig:acbd0}b). Here $i,k,l,m,\mu,\nu,\kappa\ge 0$, with
$k>0$ only if $\nu\le 1$ and $m>0$ only if $\mu\le 1$.

\noindent{\bf Type $X$.} An irreducible quadrilateral $X^\kappa_{\mu\nu}$ with digons $D_i$, $D_k$, $D_l$, $D_m$ attached to its sides
$L_1$, $L_2$, $L_3$, $L_4$ (see Fig.~\ref{fig:acbd1}a). Here $i,k,l,m,\mu,\nu,\kappa\ge 0$, with
$i>0$ only if $\mu\le 1$ and $k>0$ only if $\nu\le 1$.

\noindent{\bf Type $\bar X$.} An irreducible quadrilateral $\bar X^\kappa_{\mu\nu}$ with digons $D_i$, $D_k$, $D_l$, $D_m$ attached to its sides
$L_1$, $L_2$, $L_3$, $L_4$ (see Fig.~\ref{fig:acbd1}b). Here $i,k,l,m,\mu,\nu,\kappa\ge 0$, with
$l>0$ only if $\nu\le 1$ and $m>0$ only if $\mu\le 1$.

\noindent{\bf Type $T\nabla$.} An irreducible triangle $T_0$ with a triangle $\nabla_{jik}$ attached to its base
so that digon $D_j$ has a common side with $T_0$,
and digons $D_l$ and $D_m$ attached to the sides $T_1$ (see Fig.~\ref{fig:acbd2}a).
Here $i,j,k,l,m\ge 0$.

\noindent{\bf Type $TT$.} Irreducible triangles $T_\mu$ and $T_\nu$, with a common
integer corner at $a_3$, attached to the opposite sides of a digon $D_{2\kappa}$
with digons $D_i$ and $D_k$ attached to the bases of $T_\mu$ and $T_\nu$, and digons $D_l$
and $D_m$ attached to the sides of $T_\mu$  and $T_\nu$ (see Fig.~\ref{fig:acbd2}b).
Here $i,j,k,l,m\ge 0$, with $i>0$ only if $\mu=0$ and $k>0$ only if $\nu=0$.

\noindent{\bf Type $\overline{TT}$.} Irreducible triangles $\bar T_\mu$ and $\bar T_\nu$,
 with a common integer corner at $a_1$, attached to the opposite sides of a digon $D_{2\kappa}$
with digons $D_m$ and $D_l$ attached to the bases of $\bar T_\mu$ and $\bar T_\nu$, and digons $D_i$
and $D_k$ attached to the sides of $\bar T_\mu$  and $\bar T_\nu$ (see Fig.~\ref{fig:acbd2}c).
Here $i,j,k,l,m\ge 0$, with $i>0$ only if $\mu=0$ and $k>0$ only if $\nu=0$.

For the type $TT$ (resp., $\overline{TT}$), the integer corner $a_1$ (resp., $a_3$) of $\Q$ is mapped to a vertex $S$
(resp., $N$) of $\P$ if $\mu$ is even (resp., odd),
while its other integer corner cannot be mapped to a vertex of $\P$.
In all other cases, both integer corners of $Q$ are mapped to non-vertex points of $\P$.
\end{prop}

\begin{figure}
\centering
\includegraphics[width=4.0in]{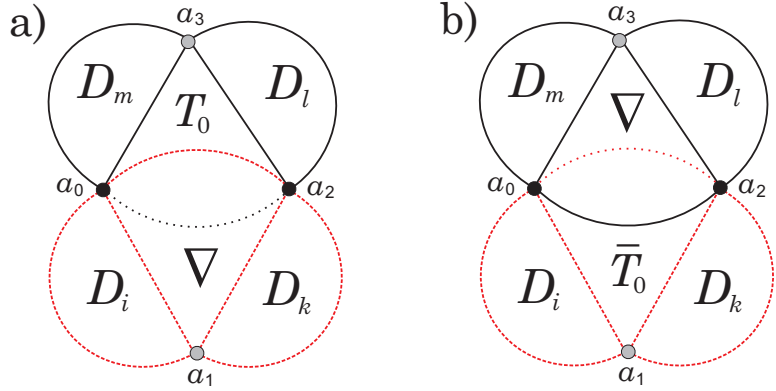}
\caption{Equivalent quadrilaterals $T\nabla$ and $\bar T\nabla$.}\label{tnabla}
\end{figure}

\begin{rmk} {\rm Note that combining $\nabla$ with $\bar T_0$ (see Fig.~\ref{tnabla}b for the $j=0$ case)
instead of $T_0$ as in type $T\nabla$ (see Fig.~\ref{tnabla}a for the $j=0$ case)
results in a quadrilateral combinatorially equivalent to the quadrilateral of type $T\nabla$.
The dotted line in Fig.~\ref{tnabla} is an interior arc of $\nabla$ which is not shown in Fig.~\ref{fig:acbd2}a.}
\end{rmk}

\section{Deformation of spherical polygons}\label{deformation}

\noindent
{\bf Deformation in a neighborhood of a corner with integer angle.}
In accordance with the previous sections, we consider the
closed upper half-plane $\H=\{ z:\Im z\geq 0\}\cup\{\infty\}$
instead of the unit disk $\D$.
Let $x\in\R=\partial H$
be a corner with an integer angle (i.e., an angle $\pi\alpha$ where $\alpha$ is integer).
and $X=f(x)$ its image under the developing map. Then the two sides adjacent
at $x$ are mapped by $f$ into a great circle, and by post-composition with a rotation
of the sphere we may assume that this great circle is the real line.
Then $f$ is real on an interval of the real line containing $x$.

Suppose that $X$ is not a vertex of $\P$.
We will show that one can deform the polygon $Q$ so that
the net does not change, and the point $X$ is
shifted to any position on some interval around $X$.

To do this, we take a disk $V$ centered at $X$ that contains no
vertices of $\P$.
Let $U$ be the component of $f^{-1}(V)$ that contains $x$.
Let $\psi_t$ be orientation preserving diffeomorphisms of $\S$ which
are equal to identity outside $V$, map $V\cap\R$ onto itself, shift $X$
to $X+t$, where $t\in\R$ is small,
and continuously depend on $t$.

Then we define
$$g_t(z)=\left\{\begin{array}{ll}\psi_t\circ f(z),& z\in U,\\
f(z),&z\in H\backslash U.\end{array}\right.$$
This is a continuous family of smooth quasiregular maps $H\to\S$, and by
the known results on solutions of Beltrami equation \cite{Ahlfors1},
one can find a continuous family of quasiconformal homeomorphisms
$\psi_t:H\to H$ such that $\phi_t\circ g_t$ are analytic.
These are the developing maps of a family of polygons which have the same
partition $\P$, and the same net $\Q$ as $Q$, but the image $X$ of the vertex $x$ of
$\Q$ is shifted on its own circle.

Notice that this procedure works also when $X$ is a vertex of $\P$,
however the net $\Q$ does change under the deformation (in a controllable way, see ...)
in the neighborhood of $x$.

In the case of a quadrilateral with two integer angles $x$ and $y$, their images $X$ and $Y$
may be on the same circle of partition $\P$,
as in Fig.~\ref{2circles-da}, when the two integer corners are adjacent, or
on two different circles, as in Fig.~\ref{2circles2a}, when the two corners are opposite.
We may assume that $S=0$ and $N=\infty$.
By a fractional linear transformation, one can fix one of the two points
$X$ or $Y$, then the second one gives a local parameter on the
set of equivalence classes of quadrilaterals. This consideration shows that
the curve in (\ref{maineq}) is non-singular for $a\not\in\{0,1\}$.
Thus it is non-singular at all real points.
\vspace{.2in}

\noindent
{\bf Degeneracy of spherical quadrilaterals}
Let us represent a marked spherical quadrilateral as a rectangle in $\C$ with vertices
$a_0,a_1,a_2,a_3$, equipped with conformal Riemannian metric with
length element
$ds=\rho(z)|dz|$ of curvature $+1$. We will call the sides $[a_0,a_1]$
and $[a_2,a_3]$ {\em horizontal} and the other two sides {\em vertical}.
Such a quadrilateral has one conformal invariant
for which we choose the extremal distance $L$ between sides $[a_1,a_2]$ and
$[a_3,a_0]$. We recall the notion of extremal length and extremal distance
\cite{Ahlfors}.

Let $\Gamma$ be
a family of curves in some region $D\subset\C$. Let $\lambda\geq 0$ be a measurable
function in $D$. We define the $\lambda$-length of a curve $\gamma$ by
$$\ell_\lambda(\gamma)=\int_\gamma\lambda(z)|dz|,$$
if the integral exists, and $\ell_\lambda(\gamma)=+\infty$ otherwise.
Then we set
$$L_\lambda(\Gamma)=\inf_{\gamma\in\Gamma}\ell_\lambda(\gamma),$$
and
$$A_\lambda(D)=\int_D\lambda^2(z)dm,$$
where $dm$ is the Euclidean area element.
Then the extremal length of $\Gamma$ is defined as
$$L(\Gamma)=\inf_\lambda\frac{L_\lambda^2(\Gamma)}{A_\lambda(D)}.$$
The extremal length is a conformal invariant. Extremal distance between two
closed sets is defined as the extremal length of the family of all curves in $D$
that connect these two sets.
For a rectangle as before, the extremal distance between
the vertical sides $[a_1,a_2]$ and
$[a_3,a_0]$ is equal to $|[a_0,a_1]|/|[a_1,a_2]|$ \cite{Ahlfors}.

In addition to the extremal distance, we consider the ordinary intrinsic distances
between the pairs of opposite sides. They are defined as the infima of $\rho$-lengths
of curves contained in our quadrilateral and connecting the two sides of a pair.

Now we have the following

\begin{lemma}\label{modulus} Consider a sequence of marked spherical quadrilaterals
whose developing maps $f$ are at most $p$-valent with a fixed integer $p$.
If the intrinsic distance between the vertical sides is bounded from below,
while the intrinsic distance between the horizontal sides tends to 0,
then the extremal distance between the vertical sides tends to $+\infty$.
\end{lemma}

{\em Proof.}
Let $\gamma_1$ be a nearly extremal curve for the intrinsic
distance between the vertical sides of $Q$,
that is the intrinsic length of $\gamma_1$ is at most $2\epsilon$.
Fix a point $P\in\gamma_1$.
Choose an arbitrary (large) number $M>0$,
and denote the intrinsic distance
by $d$.
Let $D$ be the ``annulus'' with respect to the
intrinsic metric of radii $r_1=2\epsilon$ and $r_2=M\epsilon$ centered at $P$,
that is
$$D=\{ z\in \overline{Q}: 2\epsilon\leq d(z,P)<M\epsilon\}.$$
Let $\epsilon$ be so small that $M\epsilon<\pi/2$.
 and
\begin{equation}\label{length}
2(M\epsilon+\epsilon)<c,
\end{equation}
where $c$ is a positive lower bound of the intrinsic distance between the
vertical sides.
As $f$ is $p$-valent, the intrinsic area of all intrinsic disks
$B(r)$ satisfies
\begin{equation}\label{area}
{\mathrm{area}}\ B(r)\leq 2\pi p r^2.
\end{equation}
Let $\Gamma$ be the family of curves in $Q$ connecting the horizontal sides.
Every curve $\gamma\in\Gamma$ must intersect $\gamma_1$ and both
``circles'' of the annulus $D$:
$$C_1=\{ z\in\overline{Q}: d(z,P)=2\epsilon\}\quad\mbox{and}\quad
C_2=\{ z\in\overline{Q}: d(z,P)=M\epsilon\}.$$
Thus $\gamma$ contains a curve
of the family $\Gamma'$ in $D$ which connect the inner ``circle'' $C_1$ to
the outer ``circle'' $C_2$. It follows that $L(\Gamma)\geq L(\Gamma').$
The extremal length $L(\Gamma')$ for metric annuli with metrics
satisfying (\ref{area}) is estimated in \cite[Lemma 6]{BE}:
$$L(\Gamma')\geq\frac{\log(r_2/r_1)}{32\pi}.$$
Substituting our values $r_1=2\epsilon$ and $r_2=M\epsilon$
we obtain $$a=1/L(\Gamma)<1/L(\Gamma')<(32\pi)/\log(M/2),$$
which proves the statement as $M$ is arbitrarily large.

\begin{example}\label{ex:adjacent}
{\rm Let $Q$ be one of the quadrilaterals with two adjacent non-integer corners
$a_0$ and $a_1$ (see Fig.~\ref{2circles-abcd}).

The images $X$ and $Y$ of the integer corners $a_2$ and $a_3$ of
the quadrilateral $Q$ in Fig.~\ref{2circles-abcd}c
under its developing map cannot be vertices of $\P$.
They belong to the same arc of a circle of $\P$ as shown in Fig.~\ref{2circles-da}.
The two integer corners are connected by an arc of order 1 of the net of $Q$, thus when their
images $X$ and $Y$ converge to a common point, remaining at a finite distance from the
vertices of $\P$, the intrinsic distance between the ``vertical'' sides $L_2$ and $L_4$ of $Q$
tends to 0, while the distance between its ``horizontal'' sides $L_1$ and $L_3$ does not tend to 0.
Lemma \ref{modulus} implies that the extremal distance between $L_2$ and $L_4$ tends to 0.
The corner $a_3$ of $Q$ is connected by an arc of order 1 to its corner $a_0$
but not to its corner $a_1$. Similarly, $a_2$ is connected to $a_1$ but not to $a_0$.
Lemma \ref{modulus} implies that, when either $Y$ converges to the image of $a_0$
or $X$ converges to the image of $a_1$, the extremal distance between $L_2$ and $L_4$ tends to $\infty$.

For a quadrilateral $Q$ in Fig.~\ref{2circles-abcd}a (resp., Fig.~\ref{2circles-abcd}b)
the image $Y$ of its corner $a_3$ (resp., the image $X$ of its corner $a_2$) cannot be a vertex of $\P$.
When $Y$ (resp., $X$) remains at a finite distance from vertices of $\P$, and $X$ converges to $Y$
(resp., $Y$ converges to $X$), the same argument as above implies that
the extremal distance between the sides $L_2$ and $L_4$ of $Q$ tends to 0.
If $Y$ (resp., $X$) tends to the image of $a_0$ (resp., $a_1$) and $X$ (resp., $Y$)
remains at a finite distance from the vertices of $P$,
Lemma \ref{modulus} implies that the extremal distance between the sides $L_2$ and $L_4$
of $Q$ tends to $\infty$.}
\end{example}


\begin{example}\label{ex:opposite-ab}{\rm
Let $Q$ be a quadrilateral of type $U$ (see Fig.~\ref{fig:acbd0}a)
with $\mu=\nu=0$, with opposite non-integer corners $a_0$ and $a_2$.
The images $X$ and $Y$ of its integer corners $a_1$ and $a_3$ cannot be vertices
of $\P$. The corner $a_1$ (resp., $a_3$) of $Q$ is connected by an arc of order 1 of its net
to each of its non-integer corners.
When $Y$ (resp., $X$) remains at a finite distance from the vertices of $\P$
and $X$ (resp., $Y$) converges to the image of $a_0$ (resp., $a_2$),
the distance between the sides $L_2$ and $L_4$ of $Q$ tends to $0$
while the distance between its sides $L_1$ and $L_3$ does not tend to 0.
Lemma \ref{modulus} implies that the extremal distance between $L_2$ and $L_4$ tends to 0.
When $Y$ (resp., $X$) remains at a finite distance from the vertices of $\P$
and $X$ (resp., $Y$) converges to the image of $a_2$ (resp., $a_0$),
the distance between $L_1$ and $L_3$ of $Q$ tends to $0$
while the distance between $L_2$ and $L_4$ does not tend to 0.
Lemma \ref{modulus} implies that the extremal distance between $L_2$ and $L_4$ tends to $\infty$.

Note that, for given $i,k,l,m$ and $\kappa$, this quadrilateral is combinatorially
equivalent to a quadrilateral of one of the types $\bar U,\,X,\,\bar X$ with $\mu=\nu=0$
and the same $i,k,l,m$ and $\kappa$.

Similar arguments show that, for a quadrilateral $Q$ of type $T\nabla$ (see Fig.~\ref{fig:acbd2}a),
the extremal distance between $L_2$ and $L_4$ tends to 0 when either the image $X$ of $a_1$ tends
to the image $N$ of $a_0$ or the image $Y$ of $a_3$ tends to the image $S$ of $a_2$,
and tends to $\infty$ when either $X$ tends to $S$ or $Y$ tends to $N$.}
\end{example}

\begin{example}\label{ex:opposite-aa}
{\rm Let $Q$ be a quadrilateral of type $U$ (see Fig.~\ref{fig:acbd0}a) with $\mu>0$, $\nu>0$, $i=l=0$,
with opposite non-integer corners $a_0$ and $a_2$.
Then integer corner $a_1$ (resp., $a_3$) of $Q$ is connected by an arc of order 1
of its net to $a_2$ but not to $a_0$ (resp., to $a_0$ but not to $a_2$).
The images $X$ and $Y$ of $a_1$ and $a_3$ cannot be vertices of $\P$.
When $Y$ converges to the image $N$ of $a_0$ and $X$ remains at a positive distance from
the vertices of $\P$, the distance between sides $L_1$ and $L_3$ of $Q$ tends to 0,
while the distance between $L_2$ and $L_4$ does not tend to 0.
Lemma \ref{modulus} implies that the extremal distance between $L_2$ and $L_4$ tends to $\infty$.
When $Y$ converges to the vertex $S$ of $\P$ and $X$ remains at a positive distance from
the vertices of $\P$, the distance between sides $L_1$ and $L_3$ of $Q$ tends to 0,
while the distance between $L_2$ and $L_4$ does not tend to 0,
because $a_3$ is connected by an arc of order 1 to a vertex $q\in L_1$ of the net of $Q$
that is mapped to $S$, but is not connected by an arc of order 1 to any point on $L_2$
mapped to $S$. Similarly, when $X$ converges to any of the two vertices of $\P$,
the distance between sides $L_1$ and $L_3$ of $Q$ tends to 0,
while the distance between $L_2$ and $L_4$ does not tend to 0.
Lemma \ref{modulus} implies that the extremal distance between $L_2$ and $L_4$ tends to $\infty$
for each possible degeneration of $Q$.}
\end{example}

\section{Chains of quadrilaterals}\label{chains}
Let $Q$ be a marked quadrilateral with non-integer corners $a_0$ and $a_2$,
 and with integer corners $a_1$ and $a_3$ mapped to
the points $X$ and $Y$, respectively, which are not vertices of $\P$.
Let $\Q$ be the net of $Q$.
When one of those points (say, $X$) approaches a vertex of $\P$ (say, $N$), and the combinatorial class
of $Q$ is fixed, the following options are available.
\begin{itemize}
\item[(a)] $Q$ degenerates (see section \ref{deformation}) so that the distance between its opposite sides
$L_1$ and $L_3$ tends to zero, while the distance between its sides $L_2$ and $L_4$ does not tend to 0.
Lemma \ref{modulus} implies that the extremal distance between $L_2$ and $L_4$ tends to $\infty$.
This happens when $\Q$ has an arc of order $1$ connecting $a_1$ with a
point on $L_3$ mapped to $N$,
but does not have an arc of order $1$ connecting $a_1$ to a point on $L_4$ mapped to $N$.
\item[(b)] $Q$ degenerates so that the distance between it opposite sides $L_2$ and $L_4$
tends to zero, while the distance between its sides $L_1$ and $L_3$ does not tend to 0.
Lemma \ref{modulus} implies that the extremal distance between $L_2$ and $L_4$ tends to 0.
This happens when $\Q$ has an arc of order $1$ connecting $a_1$ with a point on $L_4$
mapped to $N$,
but does not have an arc of order $1$ connecting $a_1$ to a point on $L_3$ mapped to $N$.
\item[(c)] $Q$ does not degenerate, but converges to a quadrilateral $Q'$ with the corner $a_1$ mapped to
the vertex $N$ of $\P$. This happens when $\Q$ does not have
an arc of order $1$ connecting $a_1$ with a point on one of the sides $L_3$ and $L_4$
mapped to $N$.
\item[(d)] $Q$ degenerates so that both distances, between $L_1$ and $L_3$ and between
$L_2$ and $L_4$, tend to zero. This happens when $\Q$ has arcs of order $1$ connecting
$a_1$ with points on $L_3$ and on $L_4$ mapped to $N$ (or an arc of order 1 connecting $a_1$
with $a_3$, when $a_3$ is mapped to $N$).
\end{itemize}

In the case (c), and in the corresponding cases with $Y$ instead of $X$
and/or $S$ instead of $N$, we say that $Q$ and $Q'$ are {\em adjacent}.
Note that a quadrilateral $Q'$ of type either $TT$ or $\overline{TT}$
(see Fig.~\ref{fig:acbd2}bc) has exactly one integer corner mapped to
a vertex of $\P$, and there are two distinct combinatorial classes of
quadrilaterals adjacent to $Q'$.

\begin{df}\label{df:chain}
{\rm For $k> 0$, a sequence $Q_0,Q'_1,Q_1,\dots,Q'_k,Q_k$ of
quadrilaterals with distinct combinatorial types,
where any two consecutive quadrilaterals are adjacent, and each of the terminal quadrilaterals
$Q_0$ and $Q_k$ has only one adjacent quadrilateral, is called
a {\em chain} of the {\em length} $k$.
A quadrilateral $Q_0$ having no adjacent quadrilaterals
is called a chain of length 0.

If both cases (a) and (b) are possible for degeneration of $Q_0$ and $Q_k$
then the chain is called an {\em ab-chain}.
If only the case (a) is possible, the chain is an {\em aa-chain}.
If only the case (b) is possible, the chain is a {\em bb-chain}.}
\end{df}

\begin{example}\label{ex:chain0}{\rm
It follows from Example \ref{ex:opposite-ab} that
a quadrilateral of type $U$ with $\mu=\nu=0$ is an ab-chain of length 0.
For given $i,k,l,m$ and $\kappa$, this quadrilateral is combinatorially
equivalent to a quadrilateral of one of the types $\bar U,\,X,\,\bar X$ with $\mu=\nu=0$
and the same $i,k,l,m$ and $\kappa$.

A quadrilateral of type $U$ with $i>0$ and $l>0$ is an ab-chain of length 0.
A quadrilateral of type $\bar U$ with $k>0$ and $m>0$ is an ab-chain of length 0.

A quadrilateral of type $X$ with $i>0$ and $k>0$ is an ab-chain of length 0.
A quadrilateral of type $\bar X$ with $l>0$ and $m>0$ is an ab-chain of length 0.

It follows from Example \ref{ex:opposite-aa} that
a quadrilateral of type $U$ with $\mu>0$, $\nu>0$, $i=l=0$ is an aa-chain of length 0.
Similarly, a quadrilateral of type $\bar U$ with $\mu>0$, $\nu)>0$, $k=m=0$ is a bb-chain of length 0.

It follows from Example \ref{ex:opposite-ab} that
a quadrilateral of type $T\nabla$ is an ab-chain of length 0.}
\end{example}

\begin{lemma}\label{lemma:chain0}
Any quadrilateral with opposite integer corners which is a chain of length 0
is combinatorially equivalent to one of the quadrilaterals in Example \ref{ex:chain0}.
\end{lemma}

\medskip
\noindent{\bf Proof.} One can check directly that any quadrilateral with opposite non-integer corners
other than one of the quadrilaterals in Example \ref{ex:chain0} is either
one of the quadrilaterals of types $TT$ and $\bar T\bar T$ or adjacent to a
quadrilateral of type either $TT$ or $\bar T\bar T$.

\begin{figure}
\centering
\includegraphics[width=5.0in]{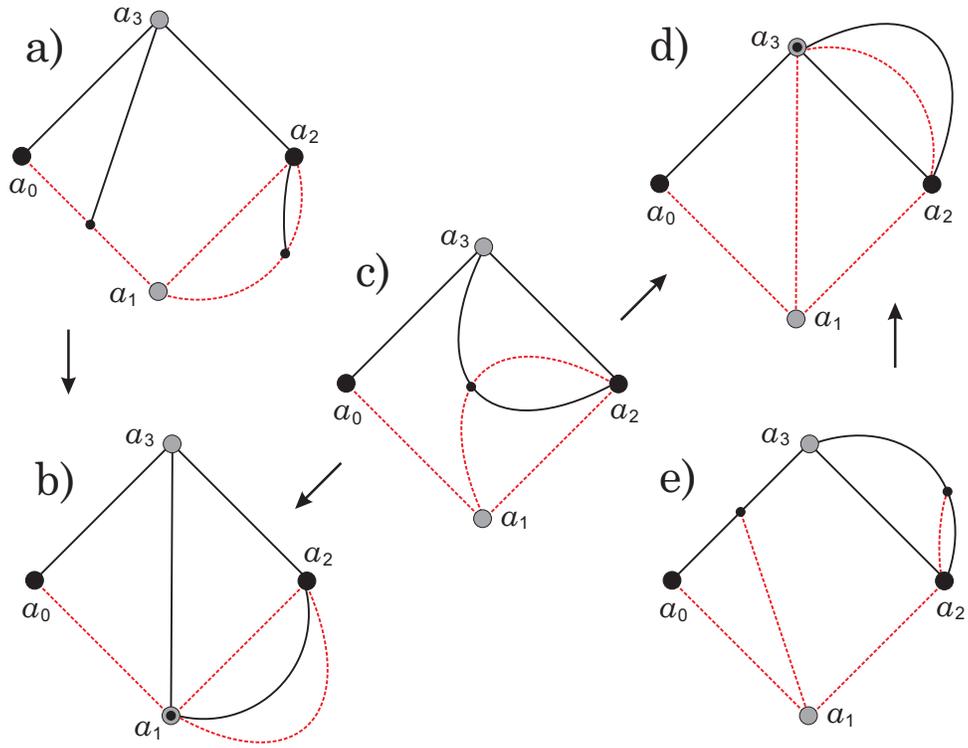}
\caption{An ab-chain of quadrilaterals of length 2.}\label{chain3}
\end{figure}

\medskip
\begin{example}\label{ex:chain3}{\rm
A chain of quadrilaterals of length $2$ is shown in Fig.~\ref{chain3}.
The quadrilateral $Q_1$ in Fig.~\ref{chain3}c is the same as the quadrilateral $I$
in Fig.~\ref{fig:non-unique}a. It can be represented either as $X_{01}$ with
a digon $D_1$ attached to its side $L_2$ or as $\bar X_{01}$ with a digon $D_1$
attached to its side $L_3$.
The quadrilateral $Q_0$ in Fig.~\ref{chain3}a is $X_{10}=U_{10}$ with a digon $D_1$ attached to
its side $L_2$.
The quadrilateral $Q_2$ in Fig.~\ref{chain3}e is $\bar X_{10}=\bar U_{10}$
with a digon $D_1$ attached to its side $L_3$.
The quadrilateral $Q'_1$ in Fig.~\ref{chain3}b is a union of two triangles $T_0$ and
a digon $D_1$.
The quadrilateral $Q'_2$ in Fig.~\ref{chain3}d is a union of two triangles $\bar T_0$ and
a digon $D_1$.

If the point $X$ to which the corner $a_1$ of $Q_0$ maps approaches $N$,
the distance between the sides $L_1$ and $L_3$ (but not of $L_2$ and $L_4$) tends to zero (case a).
The same happens if the point $Y$ to which the corner $a_3$ of $Q_0$ maps approaches $S$.
If $X$ approaches $S$, the quadrilateral $Q_0$ converges to $Q'_1$ (case c).
If $Y$ approaches $N$,
both distances (between $L_1$ and $L_3$, and between $L_2$ and $L_4$) tend to zero
(case d).
If the point $X$ to which the corner $a_1$ of $Q_1$ maps approaches $N$,
both distances (between $L_1$ and $L_3$, and between $L_2$ and $L_4$) tend to zero
(case d).
The same happens if the point $Y$ to which the corner $a_3$ of $Q_1$ maps approaches $N$.
If $X$ approaches $S$, the quadrilateral $Q_1$ converges to $Q'_1$ (case c).
If $Y$ approaches $S$, the quadrilateral $Q_1$ converges to $Q'_2$ (case c).

If the point $Y$ to which the corner $a_3$ of $Q_2$ maps approaches $N$,
the distance between the sides $L_2$ and $L_4$ (but not of $L_1$ and $L_3$) tends to zero (case b).
Thus the chain in Fig.~\ref{chain3} is an ab-chain.}
\end{example}

\begin{figure}
\centering
\includegraphics[width=5.0in]{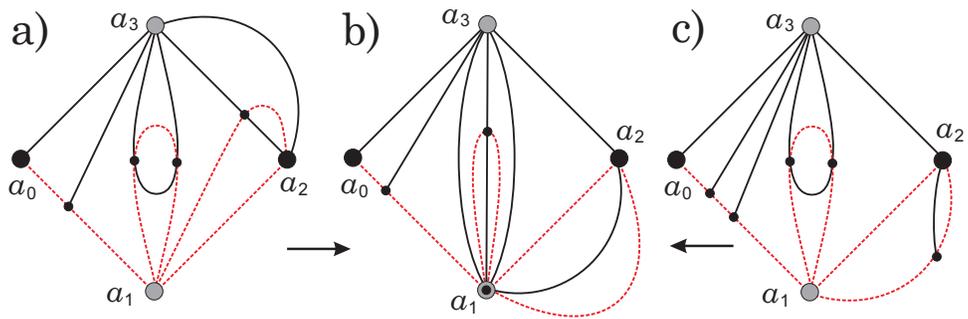}
\caption{An aa-chain of quadrilaterals of length 1.}\label{chain2}
\end{figure}

\begin{example}\label{ex:chain2}{\rm
A chain of quadrilaterals of length $1$ is shown in Fig.~\ref{chain2}.
The quadrilateral $Q'_1$ in Fig.~\ref{chain2}b has type $TT$ (see Fig.~\ref{fig:acbd2}b)
with $\mu=1$, $\nu=0$, $\kappa=1$, $i=l=m=0$ and $k=1$.

The quadrilateral $Q_0$ in Fig.~\ref{chain2}a has type $U$ (see Fig.~\ref{fig:acbd0})
with $\mu=\nu=1$, $\kappa=1$, $i=k=m=0$ and $l=1$.
When the image $X$ of the integer corner $a_1$ of $Q_0$ approaches the image $N$ of
its non-integer corner $a_0$, the quadrilateral $Q_0$ does not degenerate, and
converges to the quadrilateral $Q'_1$. When $X$ approaches the image $S$ of the non-integer
corner $a_2$ of $Q_0$, the quadrilateral $Q_0$ degenerates so that the distance
between its sides $L_1$ and $L_2$ tends to 0, while the distance between its sides
$L_2$ and $L_4$ does not tend to 0. Lemma \ref{modulus} implies
that the extremal distance between $L_2$ and $L_4$ tends to $\infty$.

The quadrilateral $Q_1$ in Fig.~\ref{chain2}c has type $U$ wit $\mu=2$, $\nu=0$, $\kappa=1$,
$i=l=m=0$, and $k=1$. It is combinatorially equivalent to a quadrilateral of type $X$ with the same values of
$\mu,\nu,\kappa,i,l,m,k$.
The same argument shows that $Q_1$ either converges to $Q'_1$ or degenerates so that
the extremal distance between its sides $L_2$ and $L_4$ tends to $\infty$.

Thus the chain in Fig.~\ref{chain2} is an aa-chain.}
\end{example}

\section{Alternative proof of Theorem \ref{theorem3}}\label{alternative}

The properties of nets described in sections \ref{intro-to-part-2} - \ref{chains}
allow us to give an alternative proof of Cases (i) and (ii) of Theorem \ref{theorem3},
computing the lower bound for the number of marked quadrilaterals with given angles and modulus,
when only two angles are non-integer.

In the case of adjacent non-integer corners (Case (i) of Theorem \ref{theorem3}),
Example \ref{ex:adjacent} implies that, moving the images of integer corners
of a quadrilateral with fixed angles in a given combinatorial equivalence class, one can
obtain a quadrilateral with the modulus (extremal distance between two of its opposite sides)
attaining any value between 0 and infinity. Thus the number of quadrilaterals with the
given angles and modulus is bounded from below by the number of distinct combinatorial classes
of the quadrilaterals.
Proposition \ref{count-abcd} implies that the number of such classes equals the number
of real solutions in Theorem \ref{theorem3} (i).

In case of opposite non-integer corners, one has to count ab-chains of quadrilaterals instead
of single combinatorial equivalence classes. It follows from Definition \ref{df:chain}
that, for any fixed angles, any admissible ab-chain $C$, and any value of the modulus,
there exists a quadrilateral combinatorially equivalent to one of the quadrilaterals in $C$
with the given angles and modulus.

Instead of counting ab-chains directly, we note that the total number of chains equals to
the number of complex solutions of equation (\ref{maineq}), which is also
the number of metrics in Theorem \ref{theorem2} and the number of
real solutions in Theorem \ref{theorem3} (i).
Indeed, due to Theorem \ref{theorem3} (iii), the number of real solutions of (\ref{maineq})
for a small $a>0$ equals the number of its complex solutions.
Every ab-chain gives one solution, and every aa-chain gives two solutions
(as it has both ends at $a=0$). Since the number of bb-chains equals the number of
aa-chains by reflection symmetry, the total number of solutions equals the total number of chains.

This implies that the number of ab-chains equals the number of all chains minus twice the number
of aa-chains, and it is enough to count aa-chains (or bb-chains).
We perform that count in this section.
Specifically, we count the aa-chains of marked quadrilaterals
having given orders $A_0,\dots,A_3$ of their corners $a_0,\dots,a_3$.

\begin{lemma}\label{lemma:chain1}
A chain of quadrilaterals is an aa-chain (resp., a bb-chain) if and only if it contains a
quadrilateral of type $U$ (resp., $\bar U$) with $\mu>0$, $\nu>0$ and $\min(i,l)=0$.
A chain may contain at most one such quadrilateral.
\end{lemma}

\medskip
\noindent{\bf Proof.}
Let us show first that each chain containing a quadrilateral $Q_0$ of type $U$ with $\mu>0$, $\nu>0$ and $\min(i,l)=0$ is an aa-chain. Example \ref{ex:opposite-aa} shows that $Q_0$ itself is
 an aa-chain of length 0 when $i=l=0$. If $i=0$ and $l>0$ then $\nu=1$ and $Q_0$ is
adjacent to a quadrilateral $Q'_1$ of type $TT$ with $\nu=0$, $l$ decreased by 1, and $k$ increased by 1.
The other quadrilateral $Q_1$ adjacent to $Q'_1$ has type $U$ with $\nu=0$, $\mu$ increased by 1,
$l$ decreased by 1, and $k$ increased by 1.
Both $Q_0$ and $Q_1$ can be degenerated so that the extremal distance between their sides $L_2$
and $L_4$ tends to $\infty$, thus $Q_0,Q'_1,Q_1$ is an aa-chain of length 1.
An example of such a chain is considered in Example \ref{ex:chain2} and shown in Fig.~\ref{chain2}.
The case when $i>0$ and $l=0$ follows by rotational symmetry.

This argument implies also that a chain containing a quadrilateral of type $U$ with $\mu>1$,
$\nu=0$, and $k>0$ (or with $\mu=0$, $\nu>1$ and $m>0$) is an aa-chain of length 1 and contains
a quadrilateral of type $U$ with $\mu>0$, $\nu>0$ and $\min(i,l)=0$.

The proof that all other chains are not aa-chains can be done case-by-case and not given here.
One of the hardest cases is considered in Example \ref{ex:chain3} and shown in Fig.~\ref{chain3}.

\begin{thm}\label{aa}
The number of aa-chains of quadrilaterals
having given orders $A_0,\dots,A_3$ of their corners $a_0,\dots,a_3$ is
\begin{equation}\label{eq:aa}
\left[\frac12\min(A_1,A_3,\delta)\right]
\end{equation}
where $\delta=\frac12\max(0,A_1+A_3-A_0-A_2)$.
\end{thm}

\medskip
\noindent{\bf Proof.}
According to Lemma \ref{lemma:chain1}, we have to count
quadrilaterals of type $U$ with $\min(\mu,\nu)>0$ and $\min(i,l)=0$
with the given orders of their corners. Note that for such a quadrilateral with $l>0$
necessarily $i=0$ and $\nu=1$, thus $\eta=(A_3-A_0)-(A_1-A_2)=\mu-\nu+2l\ge 2$.
Similarly, if $i>0$ then $l=0$ and $\mu=1$, thus $\eta=\mu-\nu-2i\le -2$.
In particular, it is enough to count
quadrilaterals with $\eta\ge 0$, for which $i=0$ and $l\ge 0$.
The case $\eta<0$ would then follow by rotation symmetry.

We start with the quadrilaterals $U_{\mu\nu}^\kappa$ with $i=l=0$, $\mu\ge\nu>0$ and
digons $D_k$ and $D_m$ attached. Then,
$$A_0=m,\;A_1=k+\nu+2\kappa+1,\;A_2=k,\;A_3=m+\mu+2\kappa+1.$$
Thus $A_0\ge 0$, $A_2\ge 0$, $A_1\ge A_2+2$, $A_3\ge A_0+2$, $A_3-A_0\ge A_1-A_2\ge 2$, and
$0\le\kappa\le \frac12(A_1-A_2-2)$.
This implies that the number of these quadrilaterals is $\left[\frac12(A_1-A_2)\right]$.

Note that each quadrilateral is uniquely determined by the value of $\kappa$,
and the possible values of $\kappa$ constitute a segment of integers
with the lower end $0$ and the upper end $\max(\kappa)$ corresponding to a quadrilateral
with $1\le\nu\le \mu$.
We'll show next that either this number equals $\left[\frac12\min(A_1,\delta)\right]$
or there exists a quadrilateral $Q'$ with the same angles as $Q$, $l>0$,
and the multiplicity of a pseudo-diagonal $\max(\kappa)+1$.

If $k=0$ then $A_2=0$ thus there are no quadrilaterals with $l>0$ and
the number of quadrilaterals is
$\left[\frac12(A_1-A_2)\right]=\left[\frac12 A_1\right]\le\left[\frac{\delta}2\right]$.

If $k=1$ then, since $A_2=1$, a quadrilateral $Q'$ with $l>0$ exists only when
$\nu=2$ and $\mu\ge 4$, in which case $\Q'$ can be taken as $U_{\mu-3,1}^{\kappa+1}$
with $D_1$ attached to the side $L_3$ and $D_m$ attached to the side $L_4$.

If $\nu=1$ then a quadrilateral $Q'$ with $l>1$ exists only when $k\ge 2$
and $\mu\ge 5$, in which case $Q'$ can be taken as $U_{\mu-4,\nu}^{\kappa+1}$
with $D_{k-2}$ attached to $L_2$, $D_2$ attached to $L_3$, and $D_m$ attached to $L_4$.

If $\nu=1$ and $\mu\le 4$ then $\left[\frac{\delta}2\right]=\kappa+1=\left[\frac12(A_1-A_2)\right]$.
If $\nu=2$ and $\mu\le 3$ then $\left[\frac{\delta}2\right]=\kappa+1=\left[\frac12(A_1-A_2)\right]$.

Next, we consider the quadrilaterals with $i=0$ and $l>0$.
For such a quadrilateral $Q$, necessarily $\nu=1$.
Let $\kappa$ be the number of pseudo-diagonals of $Q$.
Then the values $m=A_0$, $k=A_1-2\kappa-2$, $l=A_2-k=A_2-A_1+2\kappa+2$, and
$\mu=A_3-A_0-l-2\kappa-1=A_3-A_0+A_1-A_2-4\kappa-3$
are uniquely determined by $\kappa$, and the conditions $k\ge 0$, $l\ge 1$, $\mu\ge 1$
imply that $2\kappa\le A_1-2$, $2\kappa\ge A_1-A_2-1$ and $2\kappa\le\delta-2$.
Thus, for the given values of $A_0,\dots,A_3$,
the available values of $\kappa$ constitute a segment in the non-negative integers
which, if non-empty, has the upper end $\left[\frac12\min(A_1,\delta)\right]-1$.
If the lower end $\min(\kappa)$ of that segment is $0$ then $2\le A_1=k+2\le A_2+1$.
Thus there are no quadrilaterals with the same values of $A_0,\dots,A_3$ and $i=l=0$.

Otherwise, for any $0\le\kappa'<\min(\kappa)$ there is a unique quadrilateral with
the same values of $A_0,\dots,A3$ and $i=l=0$.

In any case, the total number of quadrilaterals with the given values of $A_0,\dots,A_3$
equals (\ref{eq:aa}).

\begin{figure}
\centering
\includegraphics[width=3.2in]{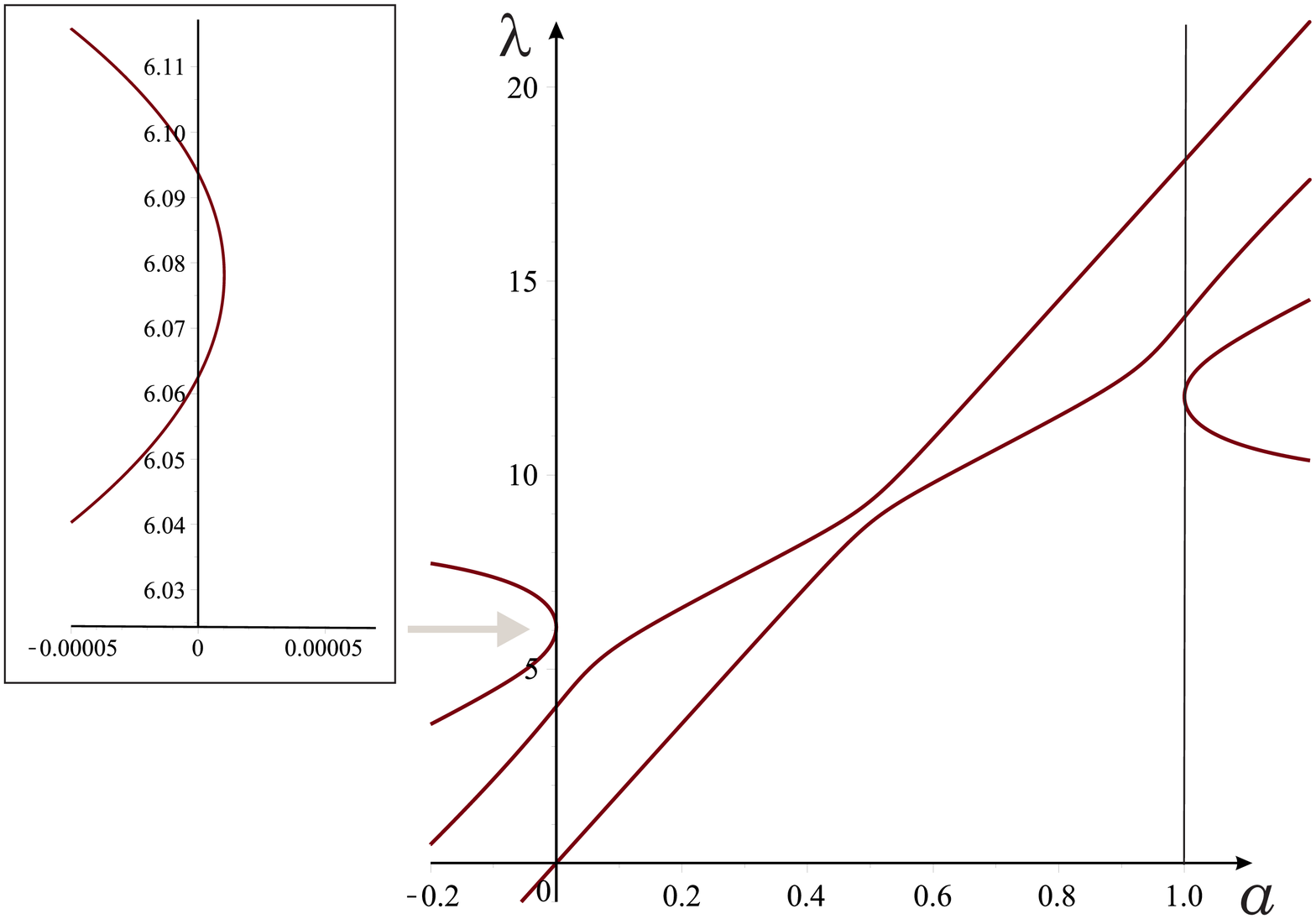}
\caption{$\alpha_1=4,\; \alpha_2=6,\; \alpha_0=\alpha_3=65/32$}\label{quad2kz}
\end{figure}

\begin{figure}
\centering
\includegraphics[width=3.2in]{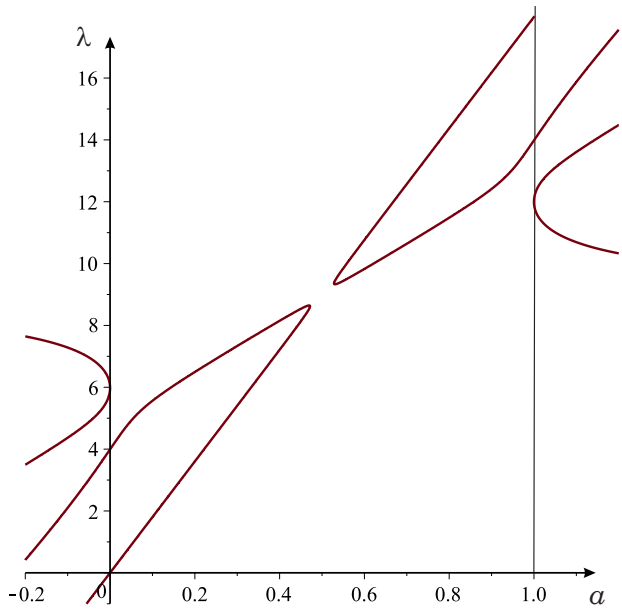}
\caption{$\alpha_1=4,\; \alpha_2=6,\; \alpha_0=\alpha_3=255/128$}\label{quad2la}
\end{figure}

\begin{figure}
\centering
\includegraphics[width=3.2in]{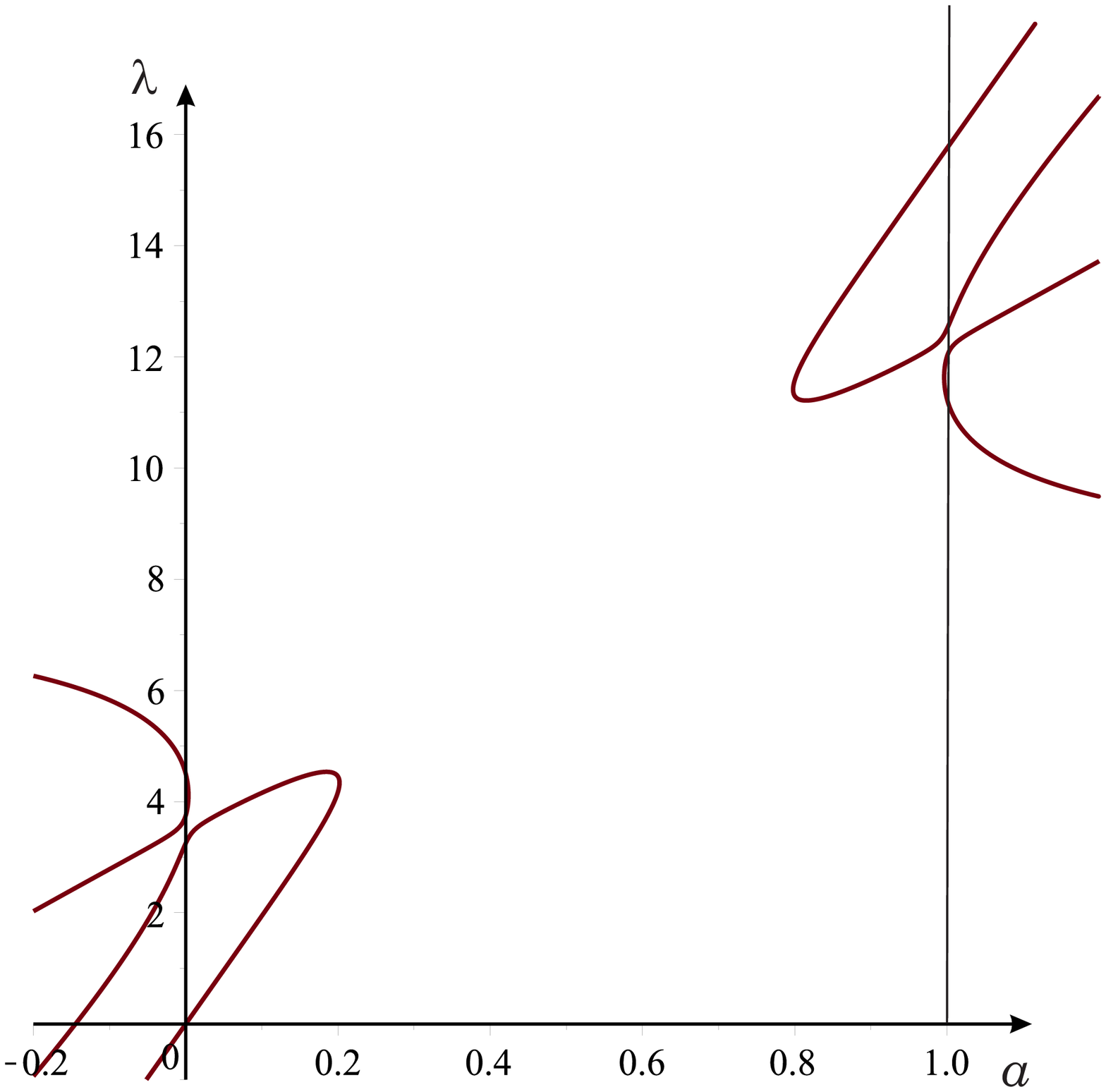}
\caption{$\alpha_1=4,\; \alpha_2=6,\; \alpha_0=\alpha_3=5/4$}\label{quad2l}
\end{figure}

\begin{figure}
\centering
\includegraphics[width=3.2in]{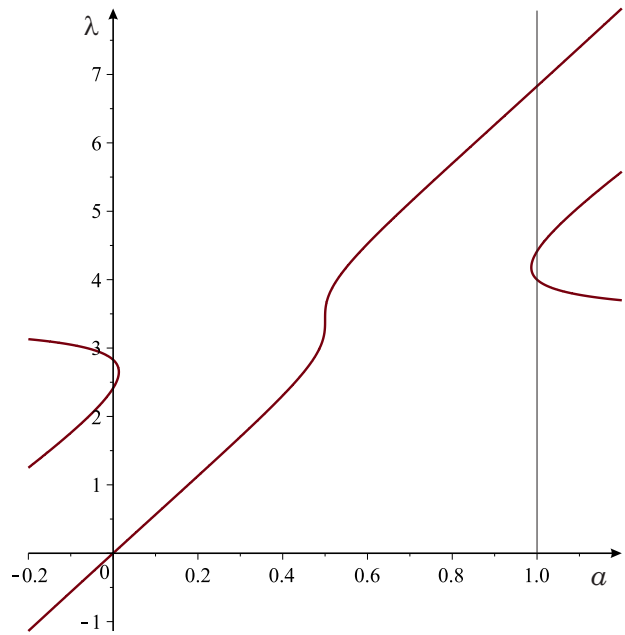}
\caption{$\alpha_1=\alpha_2=3,\; \alpha_0=\alpha_3=\sqrt{2}$}\label{quad2b}
\end{figure}

\begin{figure}
\centering
\includegraphics[width=3.2in]{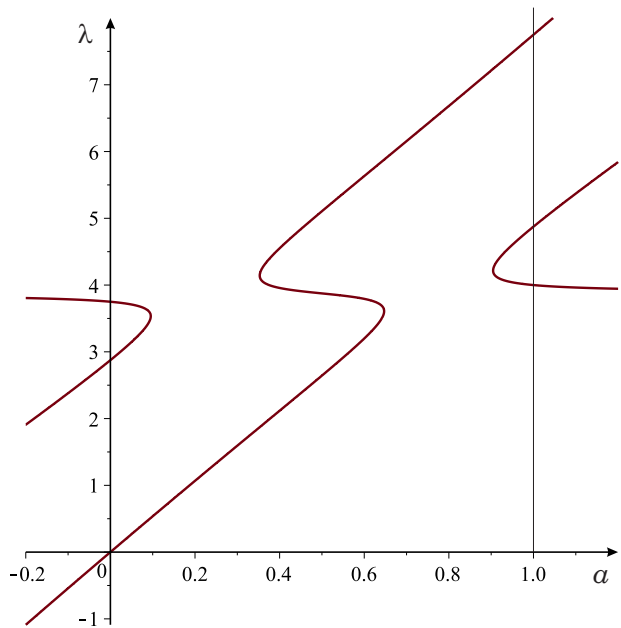}
\caption{$\alpha_1=\alpha_2=3,\; \alpha_0=\alpha_3=15/8$}\label{quad2ba}
\end{figure}

\begin{figure}
\centering
\includegraphics[width=3.2in]{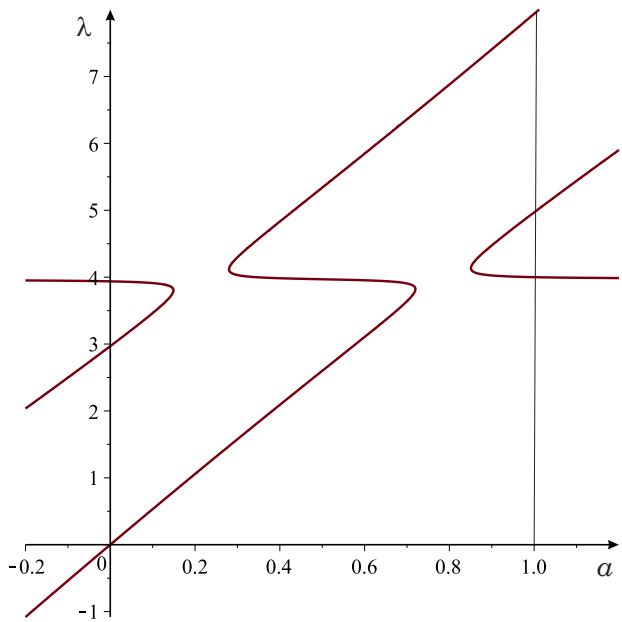}
\caption{$\alpha_1=\alpha_2=3,\; \alpha_0=\alpha_3=63/32$}\label{quad2bb}
\end{figure}

\section{Examples}\label{sec:examples}

In the following examples we choose the upper half-plane conformal model
with corners $0,1,a,\infty$, integer angles $\alpha_1$ and $\alpha_2$
at $0$ and $1$, non-integer angles $\alpha_0$ and $\alpha_3$
at $a$ and $\infty$, so the Heun equation has the form
$$y''+\left(\frac{1-\alpha_1}{z}+\frac{1-\alpha_2}{z-1}+
\frac{1-\alpha_0}{z-a}\right)y'+\frac{\alpha'\alpha''
z-\lambda}{z(z-1)(z-a)}y=0.$$
We plot the real part of the curve 
\begin{equation}\label{curva}
F(a,\lambda)=0
\end{equation}
as in (\ref{maineq}) which is defined
by the condition that the monodromy of Heun's equation is unitary.
In our examples (\ref{curva}) is the condition of absence of
logarithms in the expansion at $0$.

The values $0<a<1$ correspond to quadrilaterals with opposite integer corners.
In section \ref{intro-to-part-2}
we did show that this curve has no real singularities
when $a\not\in\{0,1\}$. That it has no singularities over $a=0$ and $a=1$
follows from the form of the Jacobi matrix: when $a=0$,
the matrix becomes triangular, with distinct diagonal entries.
In Fig.~\ref{quad2kz}, the equation (\ref{maineq})
is of degree $4$ in $\lambda$,
and it has at least $2$ real solutions for all $a$. For $a$ close to
$0$ or $1$ it has $4$ distinct real solutions.
In Figs.~\ref{quad2la} and \ref{quad2l}, there are no real solutions for some values of $a$.
Figs.~\ref{quad2b}, \ref{quad2ba} and \ref{quad2bb}
show that the number of real solutions can be larger
than the lower estimate given by Theorem~\ref{theorem3} even when $a$ is not close to
either 0 or 1.

\vspace{.1in}

{\em A. E. and A. G.: Department of Mathematics, Purdue University,

West Lafayette, IN 47907-2067 USA
\vspace{.1in}

V. T.: Department of Mathematics, IUPUI,

Indianapolis, IN 46202-3216 USA;

St. Petersburg branch of Steklov Mathematical Institute,}


\begin{thebibliography}{11}
\bibitem{Ahlfors1} L. Ahlfors, Lectures on quasiconformal mappings. Second edition, AMS, Providence, RI, 2006.
\bibitem{Ahlfors} L. Ahlfors, Conformal invariants. Topics in geometric
function theory, Reprint of the 1973 original, AMS Chelsea Publishing,
Providence, RI, 2010.
\bibitem{B} I. Biswas, A criterion for the existence of a parabolic stable
bundle of rank two over the projective line, Intl. J. Math., 9 (1998) 523--533.
\bibitem{BE} M. Bonk and A. Eremenko, Uniformly hyperbolic surfaces,
Indiana Univ. Math. J. 49 (2000) 61--80.
\bibitem{REU} R. Buckman and N. Schmitt, Spherical polygons and
unitarization, www.gang.umass.edu/reu/2002/gon.pdf.
\bibitem{CWX} Q. Chen, Y. Wu and B. Xu, Developing map and character $1$-form
of a singular conformal metric on constant curvature one on a
compact Riemann surface, arXiv:1302.6457v1.
\bibitem{D2} J. Dorfmeister and M. Schuster, Construction of planar
CMC 4-noids of genus g=0, JP Journal Geometry and Topology, 6 (2006) 3,
319--321.
\bibitem{D} J. Dorfmeister and J.-H. Eschenburg,
Real Fuchsian equations and constant mean curvature, Matem\'atica
Contempor\^{a}nea, 35 (2008), 1--25.
\bibitem{E} A. Eremenko, Metrics of positive curvature
with conic singularities on the sphere,
Proc. Amer. Math. Soc. 132 (2004), 3349--3355
\bibitem{EG0} A. Eremenko and A. Gabrielov, Rational functions with real
critical points and the B. and M. Shapiro conjecture in real enumerative
geometry,
Ann. Math., 155 (2002) 105-129.
\bibitem{EG} A. Eremenko and A. Gabrielov, Elementary proof of the
B. and M. Shapiro conjecture for rational functions,
in the book: Notions of positivity and the geometry of polynomials,
trends in mathematics, Springer, Basel, 2011, p. 167-178.
\bibitem{EG2} A. Eremenko and A. Gabrielov, Counterexamples to pole placement
by static output feedback, Linear Algebra and Appl., 351-352 (2002) 211--218.
\bibitem{EGT} A. Eremenko, A. Gabrielov and V. Tarasov,
Metrics with conic singularities and spherical polygons,
arXiv:1405.1738.
\bibitem{EGSV} A. Eremenko, A. Gabrielov, M. Shapiro and A. Vainshtein,
Rational functions and real Schubert calculus, Proc. AMS, 134
(2006), no. 4, 949--957.
\bibitem{For} L. Ford, Automorphic functions, NY, McGraw Hill, 1929.
\bibitem{FKKRUY} S. Fujimori, Y. Kawakami, M. Kokubu, W. Rossman,
M. Umehara and K. Yamada, CMC-1 trinoids in hyperbolic 3-space and metrics
of constant curvature one with conical singularities on the 2-sphere,
Proc. Japan Acad., 87 (2011), 144--149.
\bibitem{GK} F. Gantmakher and M. Krein, Oscillation matrices and kernels and
small vibrations of mechanical systems, AMS Chelsea Publ., Providence, RI,
2000.
\bibitem{G} L. Goldberg, Catalan numbers and branched
coverings by the Riemann sphere,
Adv. Math. 85 (1991), no. 2, 129--144.
\bibitem{Hilb1} E. Hilb, \"Uber Kleinsche Theoreme in der Theorie der
linearen Differentialgleichingen, Ann. Math., 66 (1909) 215--257.
\bibitem{Hilb2} E. Hilb, \"Uber Kleinsche Theoreme in der Theorie der
linearen Differentialgleichungen (2 Mitteilung),
Ann. Math., 68 (1910) 24--71.
\bibitem{HR} C. Hodgson, I. Rivin,
A characterization of compact convex polyhedra in hyperbolic 3-space,
Invent. Math. 111 (1993), no. 1, 77--111.
\bibitem{H} A. Hurwitz, \"Uber die Nullstellen der hypergeometrischen Funktion,
Math. Ann. 64 (1907) 517--560.
\bibitem{I} W. Ihlenburg, \"Uber die geometrischen Eigenschaften der
Kreisbogenvierecke, Nova Acta Leopoldina, 91 (1909) 1-79 and 5 pages of tables.
\bibitem{I2} W. Ihlenburg, Ueber die gestaltlichen Verg\"altnisse der
Kreisbogenvierecke, G\"ottingen Nachrichten, (1908) 225-230.
\bibitem{Klein} F. Klein, Ueber die Nullstellen der hypergeometrischen Reihe,
Math. Ann., 37 (1890) 573--590.
\bibitem{Klein2} F. Klein, Bemerkungen zur Theorie der
linearen Differentialgleichungen zweiter Ordnung,
Math. Ann., 64 (1907) 176--196.
\bibitem{LT} F. Luo and G. Tian,
Liouville equation and spherical convex polytopes,
Proc. Amer. Math. Soc. 116 (1992), no. 4, 1119--1129.
\bibitem{MO} R. McOwen,
Point singularities and conformal metrics on Riemann surfaces, Proc. Amer. Math. Soc. 103 (1988), 222--224.
\bibitem{MTV} E. Mukhin, V. Tarasov and A. Varchenko,
The B. and M. Shapiro conjecture in real algebraic geometry
and the Bethe ansatz,
Ann. of Math. (2) 170 (2009), no. 2, 863--881.
\bibitem{MTV2} E. Mukhin. E. Tarasov and A. Varchenko,
On reality property of Wronski maps,
Confluentes Math. 1 (2009), no. 2, 225--247.
\bibitem{P1} E. Picard, De l'equation $\Delta u=ke^u$ sur une surface de
Riemann ferm\'ee, J. Math. Pures Appl 9 (1893) 273--292.
\bibitem{P2} E. Picard, De l'equation $\Delta u=e^u$, J. Math Pures Appl.,
4 (1898) 313--316.
\bibitem{P} E. Picard, De l'integration de l'equation $\Delta u=e^u$
sur une surface de Reimann ferm\'ee, J. reine angew. Math., 130 (1905)
243--258.
\bibitem{P3} E. Picard, Quelques applications analytiques de la th\'eorie
des courbes et des surfaces alg\'ebriques, Gauthier-Villars, Paris, 1931.
\bibitem{Po} H. Poincar\'{e}, Les fonctions Fuchsiennes et l'\'equation
$\Delta u=e^u$, J. Math. pures et appl., 4 (1898) 137--230.
\bibitem{Pnt} L. Pontrjagin, Hermitian operators in spaces with
indefinite metric, Bull. Acad. Sci URSS (Izvestiya Akad. Nauk SSSR. Russian,
English summary),
8 (1944) 243-280.
\bibitem{R} A. Ronveaux, ed., Heun's differential equations, Oxford Univ. Press,
NY, 1995.
\bibitem{SG} H. P. de Saint-Gervais, Uniformisation des surfaces de Riemann,
ENS \'Editions, 2010.
\bibitem{S} I. Scherbak,
Rational functions with prescribed critical points,
Geom. Funct. Anal. 12 (2002), no. 6, 1365--1380.
\bibitem{Sch} A. Schoenflies,
Ueber Kreisbogendreiecke und Kreisbogenvierecke,
Math. Ann., 44 (1894) 105--124.
\bibitem{Sch2} A. Schoenflies, Ueber Kreisbogenbolygone,
Math. Ann. 42 (1893) 377--408.
\bibitem{Sm} V. I. Smirnov, The problem of inversion of a linear differential
equation of second order with four singularities, Petrograd, 1918 (Russian).
Reproduced in V. I Smirnov, Selected works, vol. 2, Analytic theory of ordinary
differential equations, St. Peterburg University, St. Peterburg, 1996.
\bibitem{Sm2} V. Smirnoff, Sur les \'equations diff\'erentialles lin\'eaires
du second ordre et la th\'eorie des fonctions automorphes, Bull. Sci. Math.,
45 (1921) 93--120, 126-135.
\bibitem{Ta}  G. Tarantello, Analytical, geometrical and
topological aspects of a class of mean field equations on surfaces,
Discrete Contin. Dyn. Syst. 28 (2010), no. 3, 931--973.
\bibitem{T} M. Troyanov, Prescribing curvature on compact surfaces with
conical singularities, Trans. Amer. Math. Soc., 324 (1991) 793--821.
\bibitem{Troy2} M. Troyanov,
Metrics of constant curvature on a sphere with two conical singularities,
Differential geometry (Pe\~{n}\'iscola, 1988), 296--06,
Lecture Notes in Math., 1410, Springer, Berlin, 1989.
\bibitem{UY} M. Umehara and K. Yamada, Metrics of constant curvature 1
with three conical singularities on the 2-sphere,
Illinois J. Math., 44, 1 (2000) 72--94.
\bibitem{VV} E. Van Vleck, A determination of the number of real and
imaginary roots of the hypergeometric series, Trans. Amer. Math. Soc., 3 (1902)
110-131.
\bibitem{Y} M. Yoshida, A naive-topological study of the
contiguity relations for hypergeometric function, Banach Center Publ.,
PAN, Warsawa 2005, 257--268.
\end{thebibliography}
\end{document}